\documentclass[3p,times]{elsarticle}
\usepackage{mathtools}
\usepackage{amssymb}
\usepackage{amsfonts}
\usepackage{dsfont}
\usepackage{amsthm}
\usepackage{cases}
\usepackage{bm}
\usepackage{geometry}
\usepackage{indentfirst}
\usepackage{tabularx}
\usepackage{booktabs}
\usepackage{graphicx}
\usepackage{ulem}
\usepackage{xspace}
\usepackage{xcolor}
\usepackage{multirow}
\newcommand{\bea}{\begin{eqnarray}}
\newcommand{\eea}{\end{eqnarray}}

\newcommand{\bp}{\begin{pmatrix}}
\newcommand{\ep}{\end{pmatrix}}

\newcommand{\R}{\mathds{R}}
\newcommand{\dt}{{\Delta t}}
\newcommand{\dx}{{\Delta x}}
\newcommand{\dy}{{\Delta y}}
\newcommand{\dz}{{\Delta z}}

\newcommand{\half}{\frac{1}{2}}

\newcommand{\xxi}{\boldsymbol{\xi}}
\newcommand{\mm}{\boldsymbol{m}}

\newcommand{\Aop}{\mathbb{A}} % average operator
\newcommand{\Gop}{\mathbb{G}} % implicit DG corner gradient
\newcommand{\Dop}{\mathbb{D}} % DG divergence operator
\newcommand{\Cop}{\mathbb{C}} % DG curl operator

\def\be{\begin{equation}}
\def\ee{\end{equation}}
\newcommand{\uu}{\bm{u}}
\newcommand{\oo}{\omega}
\newcommand{\boo}{\bm{\omega}}
\newcommand{\bpsi}{\bm{\Psi}}

\newcommand{\qq}{\bm{q}}
\newcommand{\xx}{\bm{x}}

\newcommand{\ijk}{_{i,j,k}}

\newcommand{\rrot}{\nabla \times \nabla \times}
\newcommand{\rot}{\nabla \times}
\newcommand{\rroth}{\nabla_h \times \nabla_h \times}
\newcommand{\roth}{\nabla_h \times}
\newcommand{\f}{\mathbf{f}}
\newcommand{\n}{\mathbf{n}}
\newcommand{\e}{\mathbf{e}}
\newcommand{\F}{\mathbf{F}}
\newcommand{\x}{\mathbf{x}}

\newcommand{\ib}{\mathbf{i}}
\newcommand{\jb}{\mathbf{j}}
\newcommand{\pb}{\mathbf{p}}

\newtheorem{theorem}{Theorem}

\newfont{\numerikEleven}{ecrm1000}
\newfont{\numerikTen}{cmss10}
\newfont{\numerikNine}{cmss9}
\newfont{\numerikEight}{cmss8}
\newfont{\numerikSeven}{cmss7}

\DeclareMathOperator{\sech}{sech}

\journal{Journal of Computational Physics}
\begin{document}
	
%!=========================================================================
%!
%!      F R O N T    M A T T E R
%!
\begin{frontmatter}
%-------------------------------------------------------
% TITLE
\title{Locally Structure-Preserving div-curl operators for high order Discontinuous Galerkin schemes}
%-------------------------------------------------------
%-------------------------------------------------------
% AUTHORS
\author[ferrara]{Walter Boscheri$^*$}
\ead{walter.boscheri@unife.it}
\cortext[cor1]{Corresponding author}

\author[ferrara]{Giacomo Dimarco}
\ead{giacomo.dimarco@unife.it}

\author[ferrara]{Lorenzo Pareschi}
\ead{lorenzo.pareschi@unife.it}

%-------------------------------------------------------
% INSTITUTIONS
\address[ferrara]{	Department of Mathematics and Computer Science \& \\ Center for Computing, Modeling and Statistics (CMCS), University of
	Ferrara, Ferrara, Italy }

%-------------------------------------------------------

%-------------------------------------------------------
% ABSTRACT 
\begin{abstract}
We propose a novel Structure-Preserving Discontinuous Galerkin (SPDG) operator that recovers at the discrete level the algebraic property related to the divergence of the curl of a vector field, which is typically referred to as div-curl problem. A staggered Cartesian grid is adopted in 3D, where the vector field is naturally defined at the corners of the control volume, while its curl is evaluated as a cell-centered quantity. Firstly, the curl operator is rewritten as the divergence of a tensor, hence allowing compatible finite difference schemes to be devised and to be proven to mimic the algebraic div-curl property. Successively, a high order DG divergence operator is built upon integration by parts, so that the structure-preserving finite difference div-curl operator is exactly retrieved for first order discretizations. We further demonstrate that the novel SPDG schemes are capable of obtaining a zero div-curl identity with machine precision from second up to sixth order accuracy. In a second part, we show the applicability of these SPDG methods by solving the incompressible Navier-Stokes equations written in vortex-stream formulation. This hyperbolic system deals with divergence-free involutions related to the velocity and vorticity field as well as to the stream function, thus it provides an ideal setting for the validation of the novel schemes. A compatible discretization of the numerical viscosity is also proposed in order to maintain the structure-preserving property of the div-curl DG operators even in the presence of artificial or physical dissipative terms. Finally, to overcome the time step restriction dictated by the viscous sub-system, Implicit-Explicit (IMEX) Runge-Kutta time stepping techniques are tailored to handle the SPDG framework. Numerical examples show that the theoretical order of convergence is reached by this new class of methods, and prove the exact preservation of the incompressibility constraint.
\end{abstract}
%-------------------------------------------------------

%-------------------------------------------------------
% KEY WORDS
%-------------------------------------------------------
\end{frontmatter}
	%!=========================================================================

%--------- SECTION --------------------------------------------------------
\section{Introduction} \label{sec:intro}
Divergence-curl (div-curl) problems are mainly originated from electrodynamics \cite{Jackson}. In particular, the Maxwell and magnetohydrodynamics (MHD) equations \cite{Davidson} involve the Faraday law, which relates the time evolution of the magnetic field to the circulation of the electric field, that is mathematically described by a curl operator. By taking the divergence of the induction equation, one gets that the divergence of the curl of the electric field obviously vanishes in order to respect the solenoidal property of the magnetic field. Another example is given by nonlinear hyperelasticity equations in solid mechanics \cite{Romenski98}, where the curl of some quantities of interest such as the distortion field must remain zero for all times if it was zero initially. These types of linear differential constraints that the governing equations exhibit at the continuous level are often labeled as involutions in literature.

In order to design numerical methods for hyperbolic conservation laws endowed with involutions, there is a need of structure-preserving schemes \cite{Pavlov}, that are able to retrieve at the discrete level the continuous properties of the governing system. Due to the important physical problems which typically deal with divergence- or curl-free conditions, an increasing interest in structure-preserving methods has emerged in the recent past. In \cite{HymanShashkov1997,JeltschTorrilhon2006,Torrilhon2004,Remi22} mimetic finite volume and finite difference schemes are proposed in a rather general context for the construction of structure-preserving schemes. Further work on finite difference schemes which respect the involutions can be found in \cite{Margolin2000,Lipnikov2014,Carney2013}. Divergence-free algorithms for plasma physics are discussed in \cite{Morrison2017,Cockburn2004}, while diffusion problems on unstructured meshes are tackled in \cite{Brezzi2005} with a mimetic finite difference approach. Because of the continuous vector space, compatible finite element discretizations are likely to be exploited for the construction of structure-preserving methods, see for instance the numerical methods presented in \cite{Nedelec1,Nedelec2,Cantarella,Hiptmair,Monk,Arnold,Alonso2015}. Focused research on the div-curl problem from a more theoretical viewpoint has been recently forwarded in \cite{Kirchhart2022}, while compatible numerical schemes have also been designed in the context of finite elements \cite{Bramble2004} and mimetic finite differences \cite{Brezzi2011}, including a weak Galerkin approach \cite{Wang2016}. A different technique is based on the use of least squares finite element methods \cite{Bensow2005,Bochev2011,Ye2020}, which generate a symmetric and positive definite system to be solved. The main advantage of these saddle point formulations is that the resulting linear system typically exhibits a
symmetric positive definite form, and thus very general and complex geometries can be easily taken into account. In this context, divergence- or curl-cleaning strategies constitute a rather different but very general approach compared to exactly structure-preserving methods, which require specific care according to the governing equations that are considered. This family of schemes was proposed in \cite{MunzCleaning,Dedneretal}, where extra terms are added to the governing system and the hyperbolic generalized Lagrangian multiplier (GLM) approach is adopted to clean those vector fields which must remain either divergence- or curl-free. Recently, this strategy has been employed in the Einstein field equations \cite{FOCCZ4GLM} or for surface tension phenomena \cite{SHTCSurfaceTension}, where it becomes extremely difficult to devise an exactly curl-free structure-preserving method. 

In the literature, structure-preserving schemes are typically concerned with the usage of staggered meshes, in order to provide natural and compatible definitions of the discrete operators such as curl, gradient and divergence. The introduction of staggered meshes was originally proposed in \cite{Yee66}, where face- and corner-staggered values of the magnetic and electric field were defined, respectively. Exactly divergence-free schemes for MHD equations have been discussed for instance in \cite{BalsaraSpicer1999,GardinerStone,SIMHD}, whereas high order structure-preserving schemes with least-squares optimization techniques have been designed in \cite{ADERdivB} for divergence-free constraints and in \cite{HOCF2021} for curl-free involutions. Corner-staggered grids have also been adopted in \cite{SIGPR} for the construction of semi-implicit structure-preserving schemes for a unified model of continuum mechanics \cite{PeshRom2014}.

The aim of this work is the construction of high order Structure-Preserving Discontinuous Galerkin (SPDG) operators that preserve the zero identity of the div-curl of any given vector field. Therefore, given a generic vector field $\f=(f_1,f_2,f_3) \in \mathds{R}^3$, with the spatial coordinate vector $\mathbf{x}=(x,y,z)=(x_1,x_2,x_3) \in \mathds{R}^3$, we want to design a discrete compatible operator such that
\begin{equation}
	\nabla_h \cdot \left( \roth \f \right) = \mathcal{O}(\epsilon),
	\label{eqn.divcurlh}
\end{equation}
where $\epsilon$ represents the machine accuracy. Contrarily to continuous finite elements, for the applications we have in mind, we prefer to work in the discontinuous Galerkin (DG) framework originally developed for neutron transport in \cite{reed} and subsequently extended to hyperbolic systems of conservation laws in \cite{cbs0,cbs1}. Therefore, the numerical representation of the unknowns is allowed to be discontinuous across cell boundaries, hence to construct a method belonging to the category of Godunov--type solvers \cite{Godunov1959} for partial differential equations (PDE). Within this setting, a structure-preserving scheme must ensure that locally within each computational cell in the mesh, the property \eqref{eqn.divcurlh} holds true in either a weak or a strong sense, or possibly in both formulations. To the best knowledge of the authors currently there exists no corner-staggered DG numerical schemes for a compatible discretization of div-curl operators in 3D. Equation \eqref{eqn.divcurlh} is computed in two steps: i) definition of a discrete curl operator, and ii) construction of a divergence operator that is compatible with the previously designed curl. At the continuous level, the curl operator explicitly writes
\begin{eqnarray}
	\left( \rot \mathbf{f} \right) &=& \left( \frac{\partial f_3}{\partial y} - \frac{\partial f_2}{\partial z}\right) \, \hat{\e}_1 + \left( \frac{\partial f_1}{\partial z} - \frac{\partial f_3}{\partial x}\right) \, \hat{\e}_2 + \left( \frac{\partial f_2}{\partial x} - \frac{\partial f_1}{\partial y}\right) \, \hat{\e}_3=\varepsilon_{\gamma \mu \tau}\hat{\e}_\tau \partial_\gamma f_\mu, \label{eqn.curl0} \\
	\left( \rot \mathbf{f} \right)_{\tau}&=& \varepsilon_{\gamma \mu \tau} \partial_\gamma f_\mu,
	\label{eqn.curl}
\end{eqnarray}
where Einstein convention implying summation over repeated indexes is assumed throughout the entire paper, and the symbol $\varepsilon_{\gamma \mu \tau}$ is the Levi-Civita tensor which can be easily computed as
\begin{equation}
	\varepsilon_{\gamma \mu \tau} = \frac{1}{2} (\gamma-\mu) \, (\mu-\tau) \, (\tau-\gamma)
\end{equation}
leading to
\begin{equation}
	\varepsilon_{123} = \varepsilon_{231}=\varepsilon_{312}=1, \quad \varepsilon_{132} = \varepsilon_{321}=\varepsilon_{213}=-1.
\end{equation}
Now, let $C$ be a control volume and $\partial C$ its boundary with an outward pointing unit normal vector $\n=(n_1,n_2,n_3) \in \mathds{R}^3$. Application of the generalized Stokes theorem to the curl definition \eqref{eqn.curl} yields
\begin{equation}
	\int \limits_{C} \varepsilon_{\gamma \mu \tau}\hat{\e}_\tau \partial_\gamma f_\mu \, dV = \int \limits_{\partial C} \varepsilon_{\gamma \mu \tau} n_\gamma f_\mu \hat{\e}_\tau \, dS. 
	\label{eqn.GST}
\end{equation}  
The cross product $\n \times \f = \varepsilon_{\gamma \mu \tau} n_\gamma f_\mu \hat{\e}_\tau$ appearing in \eqref{eqn.GST} can be also evaluated by applying Gauss theorem and defining the associated tensor $\F$:
\begin{equation}
	\int \limits_{\partial C} \varepsilon_{\gamma \mu \tau} n_\gamma f_\mu \hat{\e}_\tau \, dS =\int \limits_{\partial C} \F\cdot \n \, dS= \int \limits_{\partial C} \F_{\gamma\tau} n_{\gamma}\hat{\e}_\tau \, dS = \int \limits_{C} \partial_\gamma \F_{ \gamma\tau}\hat{\e}_\tau \, dV, \qquad \F = \left[ \begin{array}{rrr}
		0 & f_3 & -f_2 \\ -f_3 & 0 & f_1 \\ f_2 & -f_1 & 0
	\end{array} \right].
	\label{eqn.divF}
\end{equation}	
The above definition is suitable for the design of structure-preserving div-curl operators, since it only requires one single operator to be defined, that is the discrete divergence operator $\nabla_h \cdot (\cdot) = \partial_i (\cdot)_i$ acting on the control volume $C$. This definition is the starting point of the SPDG discretization discussed in this work, which will compute the discrete curl relying on corner-staggered values of the tensor field $\F$. Next, the same operator will be slightly modified to ensure the algebraic relation \eqref{eqn.divcurlh} by obtaining the discrete divergence defined at the cell centers of the control volume.

In a second part, once the SPDG div-curl operator is devised, a demonstration of applicability is proposed by considering the incompressible Navier-Stokes equations written in vortex-stream formulation. This PDE system has the characteristic of being endowed with three divergence-free conditions related to the velocity, vorticity and stream function vector fields, and the physical fluxes can be formulated by means of a rotor operator applied to a cross product. Therefore, the vortex-stream Navier-Stokes equations provide a suitable setting to properly check and validate the novel SPDG div-curl operators presented in this article. Continuous finite element methods are very popular in this field \cite{TaylorHood,SUPG,SUPG2,Fortin,Verfuerth,Rannacher1,Rannacher3} as well as finite difference schemes \cite{markerandcell,patankarspalding,patankar,vanKan}.  For what concerns the time discretization of such model, very powerful techniques based on implicit-explicit (IMEX) methods \cite{AscRuuSpi,BP2017,BosRus,PR_IMEX,BPR2017} have recently been employed to design asymptotic preserving methods for the inviscid \cite{BDLTV2020} and viscous \cite{BDT_cns,DECARIA2021113661,BosPar2021} compressible flows. Consequently, after the introduction of the SPDG method for the space derivative, we successively introduce an IMEX time stepping strategy to make the CFL--type stability condition independent of the severe restriction imposed by the viscous sub-system. Inspired by the work of Chorin \cite{Chorin67}, to further test the div-curl operators, we will propose a fully compatible discretization of artificial and physical diffusive terms, that can preserve up to machine precision the zero identity given by \eqref{eqn.divcurlh}.

The outline of this article is as follows. In Section \ref{sec:nodalDG} we introduce the discretization of the computational domain and the nodal discontinuous Galerkin framework. This setting will be used to devise the novel structure-preserving DG div-curl operators, which are fully detailed in Section \ref{sec:divcurl}. A preliminary validation test of the algebraic relation \eqref{eqn.divcurlh} is also included. Section \ref{sec:NSvorticity} is devoted to the introduction of the vortex-stream incompressible Navier-Stokes equations and the associated compatible discretization, whereas Section \ref{sec:numtest} demonstrates the applicability of the new SPDG methods to a set of test problems. A last section ends the article in which some conclusions are drawn and future investigations are proposed.

%--------- END OF SECTION -------------------------------------------------

%--------- SECTION --------------------------------------------------------
\section{Nodal discontinuous Galerkin discretization on staggered Cartesian meshes in 3D} \label{sec:nodalDG}
\begin{figure}[!htbp]
	\begin{center}
		\begin{tabular}{cc}
			\includegraphics[width=0.47\textwidth]{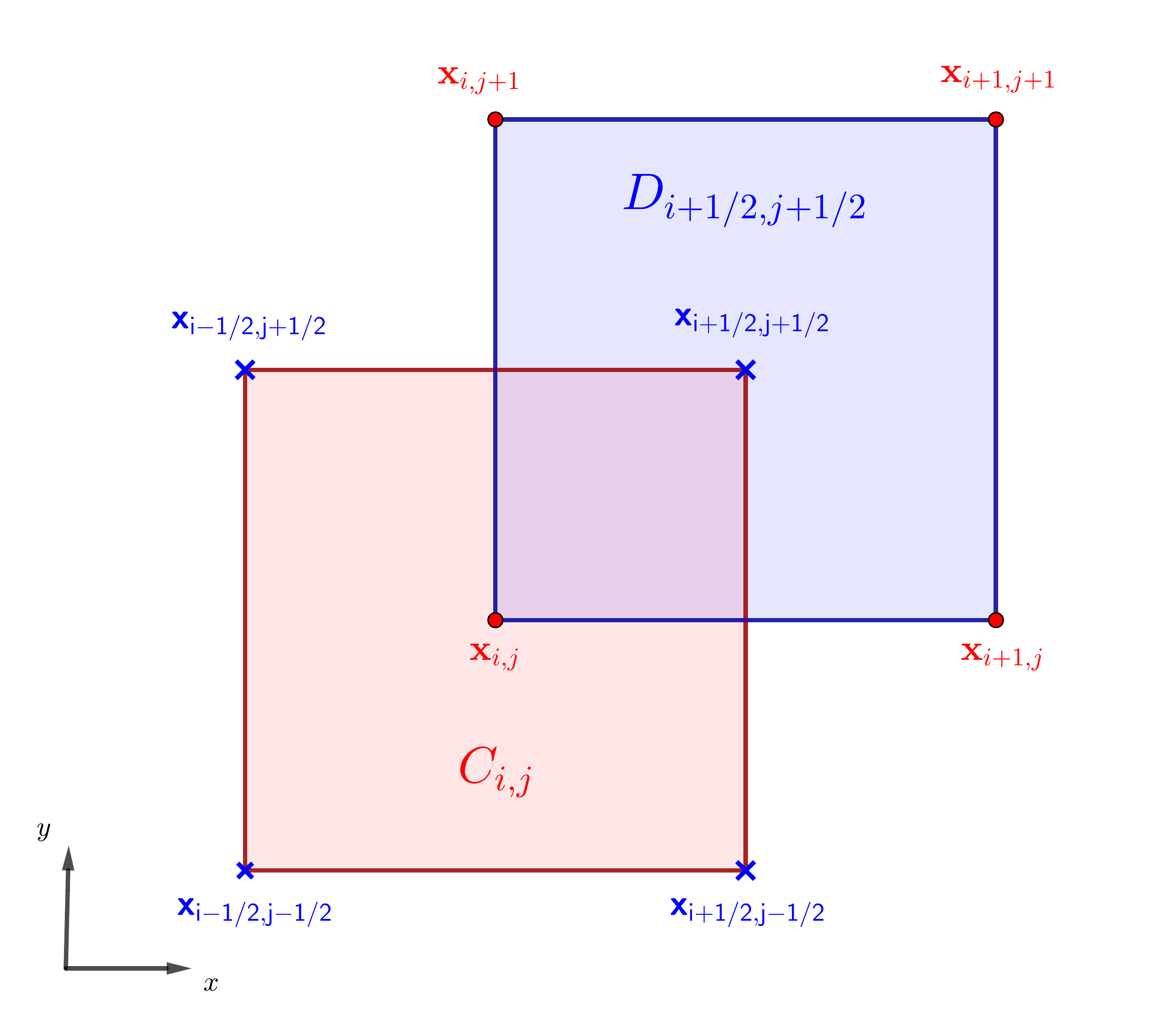}  &
			\includegraphics[width=0.47\textwidth]{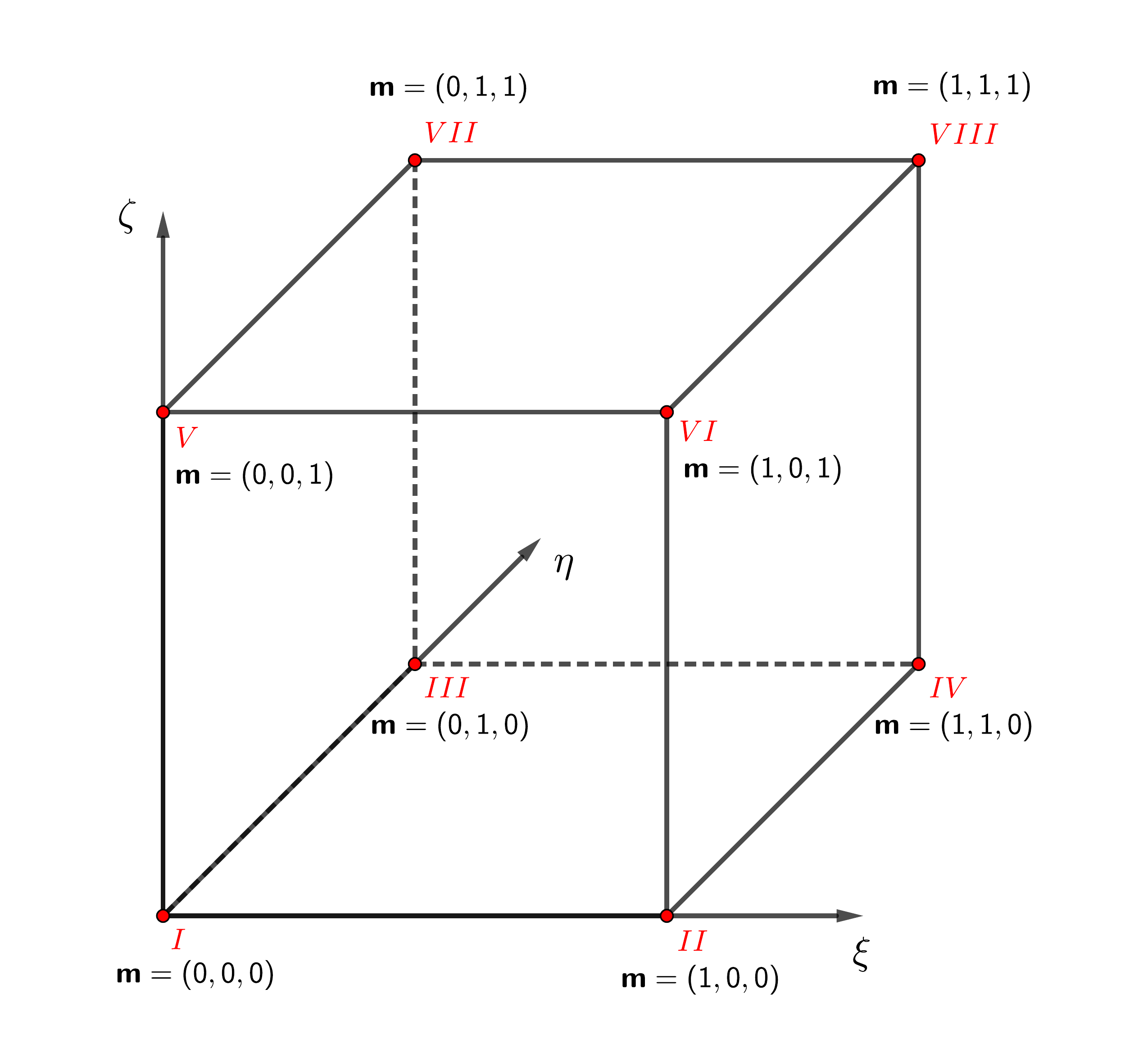} \\
		\end{tabular}
		\caption{Left: primal $C_{\ijk}$ and dual $D_{i+1/2,j+1/2,k+1/2}$ control volumes used to discretize the computational domain $\Omega(\xx)$. Right: corner orientation $\{I, \ldots, VIII \}$ and associated multi-index $\mm=(m_1,m_2,m_3)$ of the reference element $C_E=[0;1]^3$.}
		\label{fig.MeshNotation}
	\end{center}
\end{figure}

The computational domain $\Omega(\xx)=[x_{\min};x_{\max}] \times [y_{\min};y_{\max}] \times [z_{\min};z_{\max}]$ is defined in three space dimensions $d=3$, and it is bounded by the minimum and maximum value of each spatial coordinate $(\xx_{\min},\xx_{\max}) \in \mathds{R}^d$. To discretize the domain $\Omega(\xx)$, a primal Cartesian grid is employed, which is composed of a total number $N_x \times N_y \times N_z$ of cells $C\ijk$ with constant volume $|C|:=\dx \, \dy \, \dz:=\prod \limits_i \dx_i$. The characteristic mesh sizes are given by
\begin{equation}
\dx=\frac{x_{\max}-x_{\min}}{N_x}, \qquad \dy=\frac{y_{\max}-y_{\min}}{N_y}, \qquad  \dz=\frac{z_{\max}-z_{\min}}{N_z}.
\end{equation}
Let us fix some notation at the aid of Figure \ref{fig.MeshNotation} on a left in a simplified two-dimensional setting. The cell center is located at position $\xx=(x_i,y_j,z_j)$, meaning that integer indexes $(i,j,k)$ refer to cell-centered quantities. Half indexes are instead adopted for labeling quantities that are defined at midpoint locations along each spatial direction, for instance $x_{i+1/2}=x_i+\dx/2$ is the right interface in $x-$direction of cell $C_{i,j,k}$. Consequently, we make also use of a dual corner-staggered Cartesian grid, with control volumes $D_{i+1/2,j+1/2,k+1/2}$ that are centered at $\xx_{i+1/2,j+1/2,k+1/2}$, i.e. at the corners of the primal mesh. The indexes $(p,q,r)$ are used to address a generic control volume, either a primal or a dual cell. To lighten the notation, we might also use the multi-index $\ib=(i,j,k)$ and $\jb=(i+1/2,j+1/2,k+1/2)$ to label a primal and a dual cell, respectively. The same holds true for $\pb=(p,q,r)$.

Both primal and dual control volumes are hexahedra which only differ by a shifting factor of $\boldsymbol{s}=(\dx,\dy,dz)/2$, while sharing the same size, that is $|C_{\ib}|=|D_{\jb}|=|C|$. Therefore, each control volume can be mapped to the cubic reference element $C_E=[0;1]^3$ defined by the vector of reference spatial coordinates $\boldsymbol{\xi}=(\xi,\eta,\zeta)$ with $(\xi,\eta,\zeta) \in [0;1]$. The mapping $\boldsymbol{\xi}=\boldsymbol{\xi}(\x,C_{\pb})$ simply writes
\begin{equation}
\xi   = \xi(x,p) = \frac{1}{\Delta x} \left( x-x_{p-1/2} \right), \quad 
\eta  = \eta(y,q) = \frac{1}{\Delta y} \left( y-y_{q-1/2} \right), \quad 
\zeta = \zeta(z,r) = \frac{1}{\Delta z} \left( z-z_{r-1/2} \right).
\label{eqn.mapping}
\end{equation}
As shown in Figure \ref{fig.MeshNotation} on the right, the reference element counts a total number of $N_C=2^d$ corners, that are identified with a multi-index $\mm=(m_1,m_2,m_3)$ which takes into account the shifting factor with respect to the origin of the reference system. The corner orientation, that assigns a unique number $\{I, \ldots, VIII\}$ to each corner of $C_E$, and the associated multi-indexes are chosen to be as follows:
\begin{align}
I: & \quad \mm=(0,0,0), & V: & \quad \mm=(0,0,1), \nonumber \\
II: & \quad \mm=(1,0,0), & VI: & \quad \mm=(1,0,1), \nonumber \\
III: & \quad \mm=(1,1,0), & VII: & \quad \mm=(1,1,1), \nonumber \\
IV: & \quad \mm=(1,1,0), & VIII: & \quad \mm=(1,1,1), 
\label{eqn.mIdx}
\end{align}  
so that the multi-indexes also correspond to the corner reference coordinates in $\boldsymbol{\xi}$. Notice that the same multi-index can be used for defining corners and neighborhood of a dual element which are addressed with $\tilde{\mm}$, that is $\tilde{\mm}=-\mm$.

The discrete approximation of a generic scalar quantity $f(\x)$ is expressed in terms of polynomials of arbitrary degree $N$, which are written as an expansion of a set of nodal basis functions $\phi_\ell(\boldsymbol{\xi})$ with corresponding degrees of freedom $\hat{f}_\ell$:
\begin{equation}
f(\x) = \phi_\ell(\boldsymbol{\xi}) \, \hat{f}_\ell, \qquad \ell = 1, \ldots, (N+1)^3.
\label{eqn.fhat}
\end{equation}
The basis functions are obtained by tensor product in all spatial dimensions, that is
\begin{equation}
\phi_\ell(\xi,\eta,\zeta) = \phi_{l_1}(\xi) \phi_{l_2}(\eta) \phi_{l_3}(\zeta), \qquad l_1,l_2,l_3=1,\ldots,N+1, \qquad \ell=1,\ldots,(N+1)^3,
\label{eqn.psi3D}
\end{equation}
with $\ell=\ell(l_1,l_2,l_3)$ being a multi-index.The one-dimensional nodal basis $\phi_\ell(\xi)$ is given by $N+1$ linearly independent Lagrange interpolating polynomials, i.e. $\left\{\phi_\ell\right\}_{l=1}^{N+1}$, passing through a set of $N+1$ nodal points $\left\{\xi_v\right\}_{v=1}^{N+1}$, which are assumed to be the Gauss-Legendre nodes. The nodal basis exhibits by construction the following interpolation property: 
\begin{equation}
\phi_\ell(\xi_v) = \delta_{\ell v}, \qquad \ell,v=1, \ldots, N+1,
\end{equation} 
with $\delta_{\ell v}$ denoting the Kronecker delta function. Thus, the value of the approximated quantity is readily available at the nodal points, that is $f(\boldsymbol{\xi}(\x_\ell))=\hat{f}_\ell$. The discretization \eqref{eqn.fhat} is defined within each control volume, either a primal or a dual one, hence admitting discontinuities across element boundaries, which is classical in the discontinuous Galerkin framework.

Since the primal and dual cells are mapped to the same reference element according to \eqref{eqn.mapping}, the same basis functions can be used to approximate the generic quantity $f(\x)$, independently whether it is located at the cell centers or at the cell corners.

%--------- END OF SECTION -------------------------------------------------

%--------- SECTION --------------------------------------------------------
\section{Structure-Preserving div-curl operator} \label{sec:divcurl}
In order to obtain a DG operator that mimics the zero div-curl algebraic property, a discrete curl and divergence operator must be defined. However, we have seen that through the application of the generalized Stokes theorem we can reformulate the curl of a vector as the divergence of a tensor field associated to the original vector, according to Equation \eqref{eqn.divF}. As a consequence, we have only to define one single divergence operator, that will be used to compute both the curl and the divergence. To introduce the discrete operator, because of the spectral representation \eqref{eqn.fhat} given in terms of basis functions and degrees of freedom, in this section we shift the element indexes to upper index while we keep the degrees of freedom as lower indexes.

\subsection{Definition of the divergence and curl operators} \label{ssec:divcurlDG}

\begin{figure}[!htbp]
	\begin{center}
		\begin{tabular}{c}
			\includegraphics[width=0.6\textwidth]{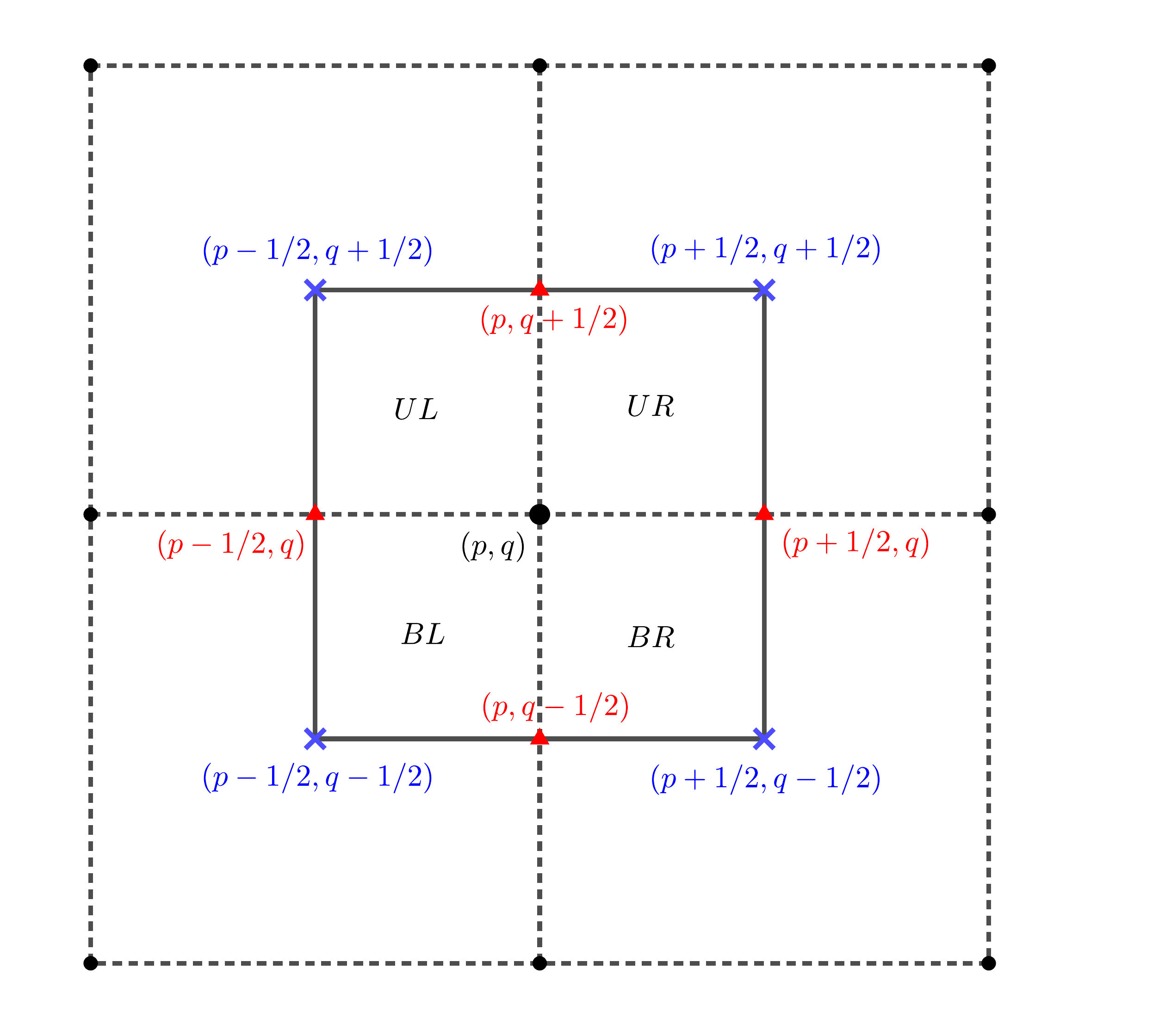}  \\
		\end{tabular}
		\caption{Construction of the divergence and gradient operators in two space dimensions, thus we set $(p,q)=(i,j)$ and $(p,q)=(i+1/2,j+1/2)$, respectively.}
		\label{fig.div2D}
	\end{center}
\end{figure}

For the sake of simplicity, let us now focus on the two-dimensional setting at the aid of Figure \ref{fig.div2D}. We consider a generic control volume $C^{p,q}$ (either a primal or a dual cell), with a vector field $\uu(\x)=(u^1,u^2) \in \mathds{R}^2$ defined at the corners, thus in the cell $C^{p+1/2,q+1/2}$. 
We aim at computing a cell-centered discrete divergence operator, that is
\begin{equation}
\boldsymbol{\alpha} = \nabla \cdot \uu, \qquad \uu(\x)=\left\{ \begin{array}{ll}
u^1(\x)&=\phi_\ell(\boldsymbol{\xi}\left(\x,C^{p+1/2,q+1/2})\right)\, \hat{u}_\ell^1 \\
u^2(\x)&=\phi_\ell(\boldsymbol{\xi}\left(\x,C^{p+1/2,q+1/2})\right)\, \hat{u}_\ell^2
\end{array}  \right. .
\label{eqn.divalpha}
\end{equation}
Adopting the DG approximation \eqref{eqn.fhat} for $\boldsymbol{\alpha}(\x)$ and using integration by parts, a weak formulation of \eqref{eqn.divalpha} over the element $C^{p,q}$ is given by
\begin{eqnarray}
\int \limits_{C^{p,q}} \phi_k \phi_\ell \, \hat{\alpha}_\ell^{p,q} \, dV &=& \int \limits_{C^{p,q}} \phi_k \, \nabla \cdot \uu \, dV \nonumber \\
&=& \int \limits_{\partial C^{p,q}} \phi_k \, \uu \cdot \n \, dS - \int \limits_{C^{p,q}} \nabla \phi_k \, \uu \, dV,
\label{eqn.divalphaW}
\end{eqnarray}
where $\phi_k$ is a test function belonging to the same polynomial space of $\phi_\ell$, and $\partial C^{p,q}$ denotes the boundary of the cell $C^{p,q}$. If $N=0$, i.e. for first order accuracy, the volume integral in \eqref{eqn.divalphaW} vanishes and only the contributions along element boundaries remain, thus retrieving a second order finite difference discretization of the divergence operator (see \ref{app.SPFD}). Due to the corner-staggered definition of the vector field $\uu$, the volume integral as well as the boundary integral in \eqref{eqn.divalphaW} must be split into all the dual elements around the same primal cell $(p,q)$, hence obtaining
\begin{eqnarray}
\int \limits_{C^{p,q}} \phi_k \phi_\ell \hat{\alpha}_\ell^{p,q} \, dV &=& \int \limits_{x_{p}}^{x_{p+1/2}} \int \limits_{y_{q}}^{y_{q+1/2}} \phi_k  \nabla \cdot \uu^{UR} dV + \int \limits_{x_{p-1/2}}^{x_{p}} \int \limits_{y_{q}}^{y_{q+1/2}} \phi_k  \nabla \cdot \uu^{UL} dV + \int \limits_{x_{p-1/2}}^{x_{p}} \int \limits_{y_{q-1/2}}^{y_{q}} \phi_k  \nabla \cdot \uu^{BL} dV + \int \limits_{x_{p}}^{x_{p+1/2}} \int \limits_{y_{q-1/2}}^{y_{q}} \phi_k  \nabla \cdot \uu^{BR} dV \nonumber \\
&+& \int \limits_{y_{q}}^{y_{q+1/2}} \phi_k \uu^{UR,\beta} \cdot \n_x^{\beta} \, dS + \int \limits_{x_{p}}^{x_{p+1/2}} \phi_k \uu^{UR,\beta} \cdot \n_y^{\beta} \, dS + \int \limits_{x_{p-1/2}}^{x_{p}} \phi_k \uu^{UL,\beta} \cdot \n_y^{\beta} \, dS - \int \limits_{y_{q}}^{y_{q+1/2}} \phi_k \uu^{UL,\beta} \cdot \n_x^{\beta} \, dS \nonumber \\ 
&-& \int \limits_{y_{q-1/2}}^{y_{q}} \phi_k \uu^{BL,\beta} \cdot \n_x^{\beta} \, dS - \int \limits_{x_{p-1/2}}^{x_{p}} \phi_k \uu^{BL,\beta} \cdot \n_y^{\beta} \, dS - \int \limits_{x_{p}}^{x_{p+1/2}} \phi_k \uu^{BR,\beta} \cdot \n_y^{\beta} \, dS + \int \limits_{y_{q-1/2}}^{y_{q}} \phi_k \uu^{BR,\beta} \cdot \n_x^{\beta} \, dS.
\label{eqn.divalphaW2}
\end{eqnarray}
The standard normals are defined as $\n_x=(1,0)^\top$ and $\n_y=(0,1)^\top$, while the index $\beta$ refers to the component of the vector $\uu$ and $\n_x$ (or $\n_y$), hence implying the contraction $\uu^{\beta}\n_x^{\beta}=\uu \cdot \n_x$. The superscripts of the sub-cells explicitly write
\begin{eqnarray}
UR=\left(p,q\right) +\half (1,1), &\qquad&
UL=\left(p,q \right)+\half (-1,1), \nonumber \\
BR=\left(p,q \right)+\half (1,-1), &\qquad&
BL=\left(p,q \right)+\half (-1,-1).
\label{eqn.indexes}
\end{eqnarray}
Now, the DG approximation \eqref{eqn.fhat} is used for the vector field $\uu$ in all its components, as written in \eqref{eqn.divalpha}, so that the spectral discretization $u^{\beta} = \phi_\ell \hat{u}^{\beta}_\ell$ is directly substituted into the weak formulation \eqref{eqn.divalphaW2}, which, after factorizing each $\uu$ belonging to the same sub-cell, leads to the following compact form of the discrete divergence operator in matrix notation:
\begin{equation}
\dx \dy \, \mathbf{M}_{k\ell} \, \hat{\alpha}_\ell^{p,q} = \mathbf{D}_{k\ell}^{UR,\beta} \hat{u}_\ell^{UR,\beta} + \mathbf{D}_{k\ell}^{UL,\beta} \hat{u}_\ell^{UL,\beta} + \mathbf{D}_{k\ell}^{BR,\beta} \hat{u}_\ell^{BR,\beta} + \mathbf{D}_{k\ell}^{BL,\beta} \hat{u}_\ell^{BL,\beta},
\label{eqn.divMatrix}
\end{equation} 
with the mass matrix given as usual by
\begin{equation}
\mathbf{M}_{k\ell} = \int\limits_{C^{p,q}}\phi_k \phi_\ell \, dV.
\label{eqn.massmatrix}
\end{equation}
Each matrix term $\mathbf{D}_{k\ell}^{*,\beta}$ in \eqref{eqn.divMatrix} contains the sub-cell volume and the two partial boundary contributions for each spatial component $\beta$ of the vector field $\uu$. To achieve a more general formulation which can ease the extension to the three-dimensional setting, let us observe that the notation used in \eqref{eqn.indexes} can be rewritten in terms of the shifting index $\mm$ introduced in \eqref{eqn.mIdx}, thus yielding
\begin{equation}
\dx \dy \, \mathbf{M}_{k\ell} \, \hat{\alpha}_\ell^{p,q} = \sum \limits_{\mm \in [-1,1]^2} \mathbf{D}_{k\ell}^{\mm,\beta} \, \hat{u}_\ell^{\mm,\beta}=\mathbf{D}_{k\ell}^{\mm,\beta} \, \hat{u}_\ell^{\mm,\beta}.
\label{eqn.divOp2D}
\end{equation} 
The vector field $\uu(\x)$ appearing in the divergence operator \eqref{eqn.divalphaW2} is defined in the physical space but in terms of the nodal basis referred to the reference space, according to \eqref{eqn.divalpha}. The test functions $\phi_k$ are defined in the reference element $C_E$, while the basis functions $\phi_\ell^{\mm}$ refer to the corner neighborhood of $C_E$, thus they must be directly evaluated as a function of the multi-index $\mm$, that is
\begin{equation}
\phi_\ell^{\mm}(\boldsymbol{\xi})=\phi_\ell \left(\boldsymbol{\xi}-\frac{\mm}{2}\right).
\end{equation}
To carry out the integration in \eqref{eqn.divalphaW2}, the intervals are mapped to the reference system by a change of variables which follows from the mapping \eqref{eqn.mapping}, hence $x=x_p+\xi \Delta x$ and $y=y_q+\eta \Delta y$. Because the integrals are evaluated numerically with Gaussian formulae of suitable accuracy \cite{stroud}, a second change of variables must be done, which shifts all the integrals to the canonical reference element defined by $\boldsymbol{\chi}=(\chi_1,\chi_2) \in [0;1]^2$, hence $\xi=(m_1+1)/4+\chi_1/2$ and $\eta=(m_2+1)/4+\chi_2/2$. The local divergence operator $\mathbf{D}_{k\ell}^{\mm,\beta}$ in \eqref{eqn.divOp2D} is then given by
\begin{eqnarray}
\mathbf{D}_{k\ell}^{m_1,m_2,\beta} &=& -\frac{1}{4} \iint_0^1 \frac{\partial }{\partial \xi^{\beta}} \phi_k \left(\frac{{\chi_1}}{2}+\frac{m_1+1}{4},\frac{{\chi_2}}{2}+\frac{m_2+1}{4}\right) \, \phi_\ell \left( \frac{1-m_2}{4} + \frac{{\chi_2}}{2}, \frac{1-m_2}{4} + \frac{\chi_2}{2}  \right) \, d{\chi_1}d{\chi_2} \, \frac{\dx \dy}{\dx^{\beta}} \nonumber \\
&& + \frac{m_1}{2} \int \limits_{0}^{1} \phi_k \left( \half + \frac{m_1}{2}, \frac{\chi_2}{2} + \frac{m_2+1}{4} \right) \, \phi_\ell \left( \half, \frac{\chi_2}{2} + \frac{1-m_2}{4} \right) \delta_{1 \beta} \, d\chi_2 \, \frac{\dx\dy}{\dx}\nonumber \\
&& + \frac{m_2}{2} \int \limits_{0}^{1} \phi_k \left( \frac{\chi_1}{2} + \frac{m_1+1}{4}, \half + \frac{m_2}{2}  \right) \, \phi_\ell \left( \frac{\chi_1}{2} + \frac{1-m_1}{4}, \half \right) \delta_{2 \beta} \,  \, d\chi_1 \, \frac{\dx\dy}{\dx},
\end{eqnarray}
with $\xi^{\beta}$ denoting the components of the reference coordinate vector $\boldsymbol{\xi}=(\xi,\eta,\zeta)=(\xi^1,\xi^2,\xi^3)$, while $(\delta_{1 \beta},\delta_{2 \beta})$ are the usual Kronecker delta functions that substitute the normal vectors $(\n_x,\n_y)$, since a Cartesian grid is employed. By means of the same philosophy (see \cite{BT22} for a similar derivation), we can extend the above operator to the general multidimensional case with $d$ space dimensions, recalling that $|C|:=\dx \dy \dz=\dx^1 \dx^2 \dx^3$:
\begin{eqnarray}
\mathbf{D}_{k\ell}^{\mm,\beta} &=& 
-\frac{|C|}{2^d\Delta x^\beta}\iiint_0^1 \frac{\partial}{\partial \xi^\beta}\phi_k\left(\frac{\mm+1}{4}+\frac{\boldsymbol{\chi}}{2}\right)\phi_\ell\left(\frac{1-\mm}{4}+\frac{\boldsymbol{\chi}}{2}\right) \, d\chi_1 d\chi_2 d\chi_3 \nonumber \\
&&+ \frac{|C| m_1}{2^{d-1}\Delta x^1}\iint_0^1
\phi_k\left(\half+\frac{m_1}{2},\frac{1+m_2}{4}+\frac{\chi_2}{2},\frac{1+m_3}{4}+\frac{\chi_3}{2}\right) 
\phi_\ell \left( \frac{1}{2}, \frac{m_2-1}{4}+\frac{\chi_2}{2},\frac{m_3-1}{4}+\frac{\chi_3}{2}\right) 
\delta_{1 \beta} d\chi_2d\chi_3 \nonumber \\
&&+ \frac{|C| m_2}{2^{d-1}\Delta x^2}\iint_0^1 
\phi_k\left(\frac{1+m_1}{4}+\frac{\chi_1}{2},\half+\frac{m_1}{2},\frac{1+m_3}{4}+\frac{\chi_3}{2}\right)  
\phi_\ell \left( \frac{m_1-1}{4}+\frac{\chi_1}{2},\frac{1}{2},\frac{m_3-1}{4}+\frac{\chi_3}{2}\right)
\delta_{2 \beta} d\chi_1d\chi_3 \nonumber \\
&&+ \frac{|C| m_3}{2^{d-1}\Delta x^3}\iint_0^1 
\phi k\left(\frac{1+m_1}{4}+\frac{\chi_1}{2},\frac{1+m_2}{4}+\frac{\chi_2}{2},\half+\frac{m_3}{2}\right) 
\phi_\ell \left(  \frac{m_1-1}{4}+\frac{\chi_1}{2},\frac{m_2-1}{4}+\frac{\chi_2}{2},\frac{1}{2}\right) 
\delta_{3 \beta} d\chi_1d\chi_2. \nonumber \\
\label{eqn.Dkl3D}
\end{eqnarray}
The divergence operator $\tilde{\Dop}_{k \ell}^{\, \pb}$ applied to a vector field with degrees of freedom $(\boldsymbol{\cdot}_{\ell}^{\mm,\beta})$ is finally derived upon multiplication by the inverse of the mass matrix and division by the cell volume, hence obtaining
\begin{equation}
\tilde{\Dop}_{k \ell}^{\, \pb}(\boldsymbol{\cdot}_{\ell}^{\mm,\beta}): = \frac{1}{\dx \dy \dz} \, \mathbf{M}_{k \ell}^{-1} \, \mathbf{D}_{k\ell}^{\mm,\beta} \, (\boldsymbol{\cdot}_{\ell}^{\mm,\beta}),
\label{eqn.divOpDG}
\end{equation}
which is referred to a cell located at spatial indexes $\pb=(p,q,r)$. The divergence operator \eqref{eqn.divOpDG} acts on the vector field $\uu(\x)$, providing the unknown degrees of freedom $\hat{\alpha}_{\ell}^{\, \pb} = \tilde{\Dop}_{k \ell}^{\, \pb} (\hat{u}^{\mm,\beta}_{\ell})$ in \eqref{eqn.divalphaW}.

According to \eqref{eqn.divF}, the discrete curl operator can be computed using the divergence operator applied to a tensor, instead of using a simple vector field. Let us define the tensor $\mathbf{U}(\x)$ associated to the vector field $\uu \in \mathds{R}^3$: 
\begin{equation}
\mathbf{U}(\x) = \left[ \begin{array}{rrr}
0 & u^3 & -u^2 \\ -u^3 & 0 & u^1 \\ u^2 & -u^1 & 0
\end{array} \right]:=\mathbf{U}^{\tau \gamma},
\end{equation}
so that, similarly to what has been done for the divergence in \eqref{eqn.divalphaW}, we can write a weak formulation of the curl:
\begin{eqnarray}
\int \limits_{C_{\pb}} \phi_k \, \varepsilon_{\gamma \mu \tau} \partial_\gamma u_\mu \, dV = \int \limits_{C_{\pb}} \phi_k \, \partial_\gamma \mathbf{U}^{\tau \gamma} \, dV.
\end{eqnarray}
The curl operator is obtained by applying the divergence operator \eqref{eqn.divOpDG} to each row of the tensor $\mathbf{U}^{\tau, \gamma}$, hence obtaining the DG curl operator $\Cop_{k\ell}^{\, \pb}(\boldsymbol{\cdot}_\ell^{\mm,\tau})$:
\begin{equation}
\Cop_{k\ell}^{\, \pb}(\boldsymbol{\cdot}_\ell^{\mm,\tau}) = \left( \roth (\boldsymbol{\cdot}_\ell^{\mm,\gamma}) \right)_{\tau}:= \tilde{\Dop}_{k\ell}^{\, \pb, \mm,\gamma}(\boldsymbol{\cdot}_{\ell}^{\mm,\tau \gamma}). 
\label{eqn.curlOpDG}
\end{equation}
The tensor is also expressed relying on the spectral approximation \eqref{eqn.fhat}, thus $\mathbf{U}(\x)=\phi_\ell(\xxi(\x,C^{\, \pb})) \, \hat{U}_\ell^{\tau \gamma}$. Likewise, the curl of the vector field $\uu(\x)$ is approximated as
\begin{equation}
	\left( \rot \uu \right)_{\tau} = \phi_\ell \, \hat{\sigma}_\ell^{\tau}. 
\end{equation}
Therefore the curl operator \eqref{eqn.curlOpDG} acts on the tensor $\mathbf{U}(\x)$ associated to a generic vector field $\uu(\x)$, providing the unknown degrees of freedom $\hat{\sigma}_\ell^{\, \pb,\tau} = \Cop_{k\ell}^{\, \pb}(\hat{\uu}_\ell^{\mm,\tau})$.

Notice that both the divergence \eqref{eqn.divOpDG} and the curl \eqref{eqn.curlOpDG} discrete operators apply to a generic cell $C^{\, \pb}$, independently of its primal or dual nature. Furthermore, each operator defines the sought degrees of freedom in the counter element, that is a cell-centered divergence $\tilde{\Dop}_{k\ell}^{\, \ib}$ starts from a corner-staggered vector field $\uu^{\, \jb}$, whereas a corner-staggered divergence $\tilde{\Dop}_{k\ell}^{\, \jb}$ takes as input a cell-centered vector field $\uu^{\, \ib}$. The same holds true for the curl operator $\Cop_{k\ell}^{\, \pb}$.  

\subsection{Structure-preserving property for DG discretization} \label{ssec:divcurlSPDG}
In \ref{app.SPFD}, we show that in the case $N=0$, i.e. in the case in which the DG discretization is reduced to a standard finite difference approximation, a structure-preserving finite difference operator is retrieved with the discretization defined in the previous section. However, it can be also shown, at the price of some tedious algebraic manipulations which we do not report here, that the discrete div-curl operator given by \eqref{eqn.divOpDG}-\eqref{eqn.curlOpDG} is not structure-preserving when $N>0$, that is
\begin{equation}
	\tilde{\Dop}_{k\ell}^{\, \ib} \left( \Cop_{k \ell}^{\, \jb} (\hat{\uu}_{\ell}^{\mm,\tau}) \right) = \kappa \, (\mathcal{O}(N+1)),
	\label{eqn.divcurlh2}
\end{equation}
meaning that the discrete condition \eqref{eqn.divcurlh} can not be satisfied up to machine accuracy $\mathcal{O}(\epsilon)$, but only up to the accuracy $\mathcal{O}(N+1)$ of the underlying DG scheme multiplied by a suitable constant $\kappa$. This is essentially due to the fact that the Schwarz theorem on the mixed second order derivative is not satisfied at the discrete level as explained later. Consequently, since we are looking for a structure-preserving divergence operator, we design here a new discrete divergence operator $\Dop_{k \ell}^{\, \pb}$ working for generic $N$, which differs from the non compatible operator $\tilde{\Dop}_{k \ell}^{\, \pb}$ (denoted with the tilde symbol) of equation \eqref{eqn.divOpDG} for a correction term.

To that aim, let us before introduce a corner gradient operator, which evaluates the gradient $\mathbf{g}=(g^1,g^2,g^3)\in \mathds{R}^3$ of a scalar field $f \in \mathds{R}$:
\begin{equation}
	\mathbf{g}=\nabla f, \qquad f = \phi_\ell(\xxi(\x,C^{\, \pb})).
\end{equation} 
Because of the duality between primal and dual mesh on a Cartesian grid, looking again at Figure \ref{fig.div2D} and following the same procedure carried out for the divergence operator, we obtain the following definition of the two-dimensional DG corner gradient operator:
\begin{equation}
	\dx \dy \, \mathbf{M}_{k\ell} \, \hat{g}_\ell^{p,q,\beta} = \sum \limits_{\mm \in \{-1,1\}^2} \mathbf{G}_{k\ell}^{\mm,\beta} \, \hat{f}_\ell^{\mm}=\mathbf{G}_{k\ell}^{\mm,\beta} \, \hat{f}_\ell^{\mm},
	\label{eqn.gradOp2D}
\end{equation} 
with the gradient matrices $\mathbf{G}_{k\ell}^{\mm,\beta}$ given by
\begin{eqnarray}
	\mathbf{G}_{k\ell}^{m_1,m_2,\beta} &=&  \frac{1}{4}\iint_0^1 \phi_k\left(\frac{m_1+1}{4}+\frac{\chi_1}{2},\frac{m_2+1}{4}+\frac{\chi_2}{2}\right)\frac{\partial}{\partial \xi^\beta}\phi_\ell\left(\frac{1-m_1}{4}+\frac{\chi_1}{2},\frac{1-m_2}{4}+\frac{\chi_2}{2}\right) \, d\chi_1 d\chi_2 \frac{\Delta x \Delta y}{\Delta x^\beta} \nonumber \\
	&& + \frac{m_1}{2}\int_0^1 \phi_k \left( \frac{1}{2}, \frac{m_2+1}{4}+\frac{\chi_2}{2}\right)\phi_\ell \left(\half-\frac{m_1}{2},\frac{1-m_2}{4}+\frac{\chi_2}{2}\right) \delta_{1 \beta} \, \, d\chi_2 \, \frac{\dx \dy}{\dx} \nonumber \\
	&& + \frac{m_2}{2}\int_0^1 \phi_k \left(\frac{m_1+1}{4}+\frac{\chi_1}{2},\half\right)\phi_ \ell\left(\frac{1-m_1}{4}+\frac{\chi_1}{2},\half-\frac{m_2}{2}\right) \delta_{2 \beta} \, \, d\chi_1 \, \frac{\dx \dy}{\dx}.
\end{eqnarray}
The extension to the three-dimensional case of the discrete gradient operator leads to
\begin{eqnarray}
	\mathbf{G}_{k\ell}^{\mm,\beta} &=& 
	\frac{|C|}{2^d\Delta x^\beta}\iiint_0^1 \frac{\partial}{\partial \xi^\beta}\phi_k\left( \frac{1-\mm}{4}+\frac{\boldsymbol{\chi}}{2} \right)\phi_\ell\left(\frac{\mm+1}{4}+\frac{\boldsymbol{\chi}}{2}\right) \, d\chi_1 d\chi_2 d\chi_3 \nonumber \\
	&&+ \frac{|C| m_1}{2^{d-1}\Delta x^1}\iint_0^1
	\phi_k\left( \frac{1}{2}, \frac{m_2-1}{4}+\frac{\chi_2}{2},\frac{m_3-1}{4}+\frac{\chi_3}{2} \right) 
	\phi_\ell \left( \half+\frac{m_1}{2},\frac{1+m_2}{4}+\frac{\chi_2}{2},\frac{1+m_3}{4}+\frac{\chi_3}{2}\right) 
	\delta_{1 \beta} d\chi_2d\chi_3 \nonumber \\
	&&+ \frac{|C| m_2}{2^{d-1}\Delta x^2}\iint_0^1 
	\phi_k\left(   \frac{m_1-1}{4}+\frac{\chi_1}{2},\frac{1}{2},\frac{m_3-1}{4}+\frac{\chi_3}{2} \right)  
	\phi_\ell \left( \frac{1+m_1}{4}+\frac{\chi_1}{2},\half+\frac{m_1}{2},\frac{1+m_3}{4}+\frac{\chi_3}{2}\right)
	\delta_{2 \beta} d\chi_1d\chi_3 \nonumber \\
	&&+ \frac{|C| m_3}{2^{d-1}\Delta x^3}\iint_0^1 
	\phi k\left(  \frac{m_1-1}{4}+\frac{\chi_1}{2},\frac{m_2-1}{4}+\frac{\chi_2}{2},\frac{1}{2} \right) 
	\phi_\ell \left( \frac{1+m_1}{4}+\frac{\chi_1}{2},\frac{1+m_2}{4}+\frac{\chi_2}{2},\half+\frac{m_3}{2}\right) 
	\delta_{3 \beta} d\chi_1d\chi_2 \nonumber \\
	&:=& -\mathbf{D}_{\ell k}^{-\mm,\beta},
	\label{eqn.Gkl3D}
\end{eqnarray}
thus the corner gradient operator can be easily computed from the divergence operator \eqref{eqn.Dkl3D}, and it compactly writes
\begin{equation}
	{\Gop}_{k \ell}^{\, \pb} (\boldsymbol{\cdot}_{\ell}^{\mm})=\nabla_h (\boldsymbol{\cdot}_{\ell}^{\mm}): = \frac{1}{\dx \dy \dz} \, \mathbf{M}_{k \ell}^{-1} \, \mathbf{G}_{k\ell}^{\mm,\beta}.
	\label{eqn.gradOpDG}
\end{equation}
Furthermore, this implies that both $\mathbf{D}_{k\ell}^{\mm,\beta}$ and $\mathbf{G}_{k \ell}^{\mm,\beta}$ share the same discrete computations, which will play a crucial role for retrieving the structure-preserving property of the DG div-curl operator. 

In fact, the reason why the div-curl discrete operator \eqref{eqn.divcurlh2} is not structure-preserving lies in the non-vanishing mixed second derivative terms, i.e. Schwarz theorem is not respected at the discrete level. Indeed, using the relation \eqref{eqn.Gkl3D} between gradient and divergence operators, we can exactly quantify the discrete errors arising from \eqref{eqn.divcurlh2}, that is
\begin{eqnarray}
	\kappa \, (\mathcal{O}(N+1)) &=& \left( \partial^{\alpha} \partial^{\beta} u^{\tau}_{\ell} - \partial^{\beta} \partial^{\alpha} u^{\tau}_{\ell} \right) (1-\delta_{\alpha \beta}) \nonumber \\
	&=&\left( \frac{1}{\dx \dy \dz} \, \mathbf{M}_{k \ell}^{-1} \mathbf{G}_{k\ell}^{\mm,\alpha} \left( \frac{1}{\dx \dy \dz} \, \mathbf{M}_{k \ell}^{-1}\mathbf{G}_{k\ell}^{\mm,\beta} (\hat{u}_{\ell}^{\mm,\tau}) \right) - \frac{1}{\dx \dy \dz} \, \mathbf{M}_{k \ell}^{-1} \mathbf{G}_{k\ell}^{\mm,\beta} \left( \frac{1}{\dx \dy \dz} \, \mathbf{M}_{k \ell}^{-1} \mathbf{G}_{k\ell}^{\mm,\alpha} (\hat{u}_{\ell}^{\mm,\tau}) \right) \right) (1-\delta_{\alpha \beta}) \nonumber \\
	&=& \frac{1}{\dx \dy \dz} \, \mathbf{M}_{k \ell}^{-1} \mathbf{D}_{k\ell}^{\mm,\tau} \cdot \left( \varepsilon_{\gamma \mu \tau} \frac{1}{\dx \dy \dz} \, \mathbf{M}_{k \ell}^{-1} \mathbf{G}_{k\ell}^{\mm,\gamma} \, \hat{u}_{\ell}^{\mm,\mu}\right),
	\label{eqn.corrDG}
\end{eqnarray} 
with the Hessian tensor given by $\mathbf{H}^{\alpha, \beta, \tau}_{\ell} \hat{u}^{\tau}_{\ell} = \mathbf{G}_{k\ell}^{\mm,\alpha} \left( \mathbf{G}_{k\ell}^{\mm,\beta} (\hat{u}_{\ell}^{\mm,\tau}) \right)$.
%Indeed, the result of the operator \eqref{eqn.divcurlh2} applied to a vector field $\uu=(u^1,u^2,u^3) \in \mathds{R}^3$ exactly amounts to the contribution of the discrete errors related to the following relations:
%\begin{eqnarray}
%	\kappa(\mathcal{O}(N+1)) &=& \left( \frac{\partial u^3}{\partial xy} - \frac{\partial u^3}{\partial yx} \right)_h + \left( \frac{\partial u^2}{\partial zx} - \frac{\partial u^3}{\partial xz} \right)_h + \left( \frac{\partial u^1}{\partial yz} - \frac{\partial u^1}{\partial zy} \right)_h \nonumber \\
%	&=& \left [ \frac{\partial}{\partial_x} \left( \frac{\partial u^3}{\partial y} - \frac{\partial u^2}{\partial z}\right) \right]_h+ \left[ \frac{\partial}{\partial_y} \left( \frac{\partial u^1}{\partial z} - \frac{\partial u^3}{\partial x}\right) \right]_h + \left[ \frac{\partial}{\partial_z} \left( \frac{\partial u^2}{\partial x} - \frac{\partial u^1}{\partial y}\right) \right]_h,
%\end{eqnarray}
%which is the discrete divergence of the curl of the vector field $\uu$. 
Thus, eventually, we can obtain a structure-preserving div-curl operator by defining the divergence operator $\Dop_{k \ell}^{\, \pb}$ as follows:
\begin{eqnarray}
	\Dop_{k \ell}^{\, \pb} (\boldsymbol{\cdot}_{\ell}^{\mm,\beta}) = \nabla_h \cdot (\boldsymbol{\cdot}_{\ell}^{\mm,\beta}) &:=& \frac{1}{\dx \dy \dz} \, \mathbf{M}_{k \ell}^{-1} \, \left( \mathbf{D}_{k\ell}^{\mm,\beta}(\boldsymbol{\cdot}_{\ell}^{\mm,\beta}) - \mathbf{D}_{k\ell}^{\mm,\beta} \cdot \left( \frac{1}{\dx \dy \dz} \, \mathbf{M}_{k \ell}^{-1} \varepsilon_{\gamma \mu \beta} \, \mathbf{G}_{k\ell}^{\mm,\gamma}(\boldsymbol{\cdot}_{\ell}^{\mm,\mu}) \right) \right) \nonumber \\
	&:=&\tilde{\Dop}_{k \ell}^{\, \pb}(\boldsymbol{\cdot}_{\ell}^{\mm,\beta}) - \frac{1}{\dx \dy \dz} \, \mathbf{M}_{k \ell}^{-1} \mathbf{D}_{k\ell}^{\mm,\beta} \cdot \left( \frac{1}{\dx \dy \dz} \, \mathbf{M}_{k \ell}^{-1} \varepsilon_{\gamma \mu \beta} \, \mathbf{G}_{k\ell}^{\mm,\beta}(\boldsymbol{\cdot}_{\ell}^{\mm,\mu}) \right).
	\label{eqn.divOpSPDG}
\end{eqnarray}
%At the continuous level, this simply corresponds to compute
%\begin{equation}
%	\nabla \dot (\boldsymbol{\cdot}) = \nabla \dot (\boldsymbol{\cdot}) - \nabla \cdot \rot (\boldsymbol{\cdot}),
%\end{equation}
%thus we formally add a zero term to the continuous divergence, which nevertheless allows the discrete errors related to the div-curl operator to vanish.

\begin{theorem}\label{lemma_divcurlDG}
	The discrete curl $\Cop_{h}=\left( \roth (\boldsymbol{\cdot}_\ell^{\mm,\tau}) \right)_{\beta}$ and the discrete divergence $\Dop_{h}=\nabla_h \cdot (\boldsymbol{\cdot}^{\mm,\beta}_{\ell})$ operators, defined by  \eqref{eqn.curlOpDG} and \eqref{eqn.divOpSPDG}, respectively, satisfy the discrete div-curl property \eqref{eqn.divcurlh} locally, namely 
	\begin{equation}
		\nabla_h \cdot \left( \left(\roth (\boldsymbol{\cdot}_\ell^{\mm,\tau})\right)_{\beta} \right)^{\mm,\beta}_{\ell}=\mathcal{O}(\epsilon) \qquad \forall \ell= 1, \ldots, (N+1)^3.
		\label{eqn.divcurlhL}
	\end{equation}
\end{theorem}

\begin{proof}
Using the relation between the gradient and the divergence operators \eqref{eqn.Gkl3D}, we can write the curl operator in terms of the gradient operator, namely
\begin{equation}
	=\left( \roth (\boldsymbol{\cdot}_\ell^{\mm,\tau}) \right)_{\beta} = \frac{1}{\dx \dy \dz} \, \mathbf{M}_{k \ell}^{-1} \varepsilon_{\gamma \mu \beta} \, \mathbf{G}_{k\ell}^{\mm,\gamma} (\boldsymbol{\cdot}_{\ell}^{\mm,\mu}), \qquad \forall \, k,\ell = 1, \ldots, (N+1)^3. 
\end{equation}
Applying the discrete divergence \eqref{eqn.divOpSPDG} to the above curl yields
\begin{equation}
	\frac{1}{\dx \dy \dz} \, \mathbf{M}_{k \ell}^{-1} \mathbf{D}_{k\ell}^{\mm,\beta} \left( \frac{1}{\dx \dy \dz} \, \mathbf{M}_{k \ell}^{-1} \varepsilon_{\gamma \mu \beta} \, \mathbf{G}_{k\ell}^{\mm,\gamma}(\boldsymbol{\cdot}_{\ell}^{\mm,\mu}) \right) - \frac{1}{\dx \dy \dz} \, \mathbf{M}_{\ell k}^{-1} \mathbf{D}_{k \ell}^{\mm,\beta} \left( \frac{1}{\dx \dy \dz} \, \mathbf{M}_{k \ell}^{-1} \varepsilon_{\gamma \mu \beta} \, \mathbf{G}_{k\ell}^{\mm,\gamma}(\boldsymbol{\cdot}_{\ell}^{\mm,\mu}) \right) = \mathcal{O}(\epsilon).
\end{equation}
Thus, because of the adoption of the nodal basis \eqref{eqn.psi3D} in the DG approximation \eqref{eqn.fhat}, the zero div-curl property holds true for each degree of freedom $(\hat{\boldsymbol{\cdot}})_\ell$, hence making the DG div-curl operator locally structure-preserving.
\end{proof}

\subsection{Numerical validation of the SPDG div-curl operator} \label{ssec.SPDGdivcurl}
In this part, we numerically demonstrate the compatible div-curl DG operator \eqref{eqn.divcurlhL}. Let us consider the following vector field $\boldsymbol{\Psi}(\x)$ and the corresponding curl given by $\uu=\rot \boldsymbol{\Psi}$:
\begin{eqnarray}
	\boldsymbol{\Psi}(\x) &=& \left\{ \begin{array}{lcl}
		\Psi_1 &=& -\sin(2\pi x) \, \cos(2\pi y) \, \cos(2\pi z) \\
		\Psi_2 &=& 2\cos(2\pi x) \, \sin(2\pi y) \, \cos(2\pi z) \\
		\Psi_3 &=& \phantom{-} \cos(2\pi x) \, \cos(2\pi y) \,  \cos(2\pi z) 
	\end{array} \right. , \\
	\boldsymbol{\uu}(\x) &=& \left\{ \begin{array}{lcl}
		u_1 &=& \phantom{-} 2 \pi \cos(2\pi x) \, \sin(2\pi y) \, \sin(2\pi z) \\
		u_2 &=& \phantom{-} 4 \pi \sin(2\pi x) \, \cos(2\pi y) \, \sin(2\pi z) \\
		u_3 &=&          -  6 \pi \sin(2\pi x) \, \sin(2\pi y) \,  \cos(2\pi z) 
	\end{array} \right. .
\end{eqnarray}
The computational domain is given by $\Omega=[0;1]^3$, and a sequence of refined computational meshes are used to carry out a convergence study. Each mesh is made of $N_h^3=N_x \times N_y \times N_z$ number of cells, considering four grids defined by $N_h=\{6,12,24,48\}$. The errors are measured in $L_1$ and $L_{\infty}$ norms, and they are related to the approximation of the vector field $\boldsymbol{\Psi}$ as well as its curl, i.e. the vector $\uu$. We expect to achieve a convergence rate of $\mathcal{O}(N+1)$ for $\boldsymbol{\Psi}$, whereas an order of $\mathcal{O}(N)$ is predicted for the vector field $\uu$, because it is computed as the curl of $\boldsymbol{\Psi}$, hence involving a derivative. Table \ref{tab.divcurl_conv} confirms that the formal accuracy of the structure-preserving DG schemes is obtained from second up to sixth order. Furthermore, we also monitor the $L_{\infty}$ norm of the div-curl error computed by considering all the degrees of freedom, hence ensuring that the SPDG schemes are locally structure-preserving up to machine precision as proven by Theorem \ref{lemma_divcurlDG}. 

\begin{table}[!htbp]  
	\caption{Numerical convergence results of the SPDG div-curl operators from second up to sixth order of accuracy. The divergence errors are also reported for each simulation and are labeled with $\left|\nabla_h \cdot (\roth \boldsymbol{\Psi})\right|_{\infty}$.}
	\begin{center} 
		\scalebox{0.95}{\begin{tabular}{c|cccccccc|c}
			\multicolumn{10}{c}{SPDG $\mathcal{O}(2)$} \\
			$N_h$ & $u_{1,L_1}$ & $\mathcal{O}(L_1)$ & $u_{1,L_{\infty}}$ & $\mathcal{O}(L_{\infty})$ & $\Psi_{1,L_1}$ & $\mathcal{O}(L_1)$ & $\Psi_{1,L_{\infty}}$ & $\mathcal{O}(L_{\infty})$ & $\left|\nabla_h \cdot (\roth \boldsymbol{\Psi})\right|_{\infty}$ \\
			\hline
			 6 & 3.0295E-01 &    - & 1.4230E+00 &    - & 1.4666E-02 &    - & 9.9726E-02 &    - & 2.64E-15 \\
			12 & 1.3426E-01 & 1.17 & 6.8122E-01 & 1.06 & 3.9319E-03 & 1.90 & 3.0475E-02 & 1.71 & 2.87E-15  \\
			24 & 6.3299E-02 & 1.08 & 3.2428E-01 & 1.07 & 9.7355E-04 & 2.01 & 8.3216E-03 & 1.87 & 1.03E-15  \\
			48 & 3.1019E-02 & 1.03 & 1.5999E-01 & 1.02 & 2.4278E-04 & 2.00 & 2.1263E-03 & 1.97 & 6.88E-16  \\
			\multicolumn{10}{c}{} \\
			\multicolumn{10}{c}{SPDG $\mathcal{O}(3)$} \\
			$N_h$ & $u_{1,L_1}$ & $\mathcal{O}(L_1)$ & $u_{1,L_{\infty}}$ & $\mathcal{O}(L_{\infty})$ & $\Psi_{1,L_1}$ & $\mathcal{O}(L_1)$ & $\Psi_{1,L_{\infty}}$ & $\mathcal{O}(L_{\infty})$ & $\left|\nabla_h \cdot (\roth \boldsymbol{\Psi})\right|_{\infty}$ \\
			\hline
			 6 & 4.0804E-02 &    - & 2.3403E-01 &    - & 1.8588E-03 &    - & 5.0871E-03 &    - & 2.35E-15 \\
			12 & 1.0351E-02 & 1.98 & 6.8490E-02 & 1.77 & 2.3513E-04 & 2.98 & 4.7027E-04 & 3.44 & 9.34E-16  \\
			24 & 2.5876E-03 & 2.00 & 1.7972E-02 & 1.93 & 2.9441E-05 & 3.00 & 7.4393E-05 & 2.66 & 2.05E-15  \\
			48 & 6.4702E-04 & 2.00 & 4.5482E-03 & 1.98 & 3.6823E-06 & 3.00 & 9.0690E-06 & 3.04 & 7.56E-15  \\
			\multicolumn{10}{c}{} \\
			\multicolumn{10}{c}{SPDG $\mathcal{O}(4)$} \\
			$N_h$ & $u_{1,L_1}$ & $\mathcal{O}(L_1)$ & $u_{1,L_{\infty}}$ & $\mathcal{O}(L_{\infty})$ & $\Psi_{1,L_1}$ & $\mathcal{O}(L_1)$ & $\Psi_{1,L_{\infty}}$ & $\mathcal{O}(L_{\infty})$ & $\left|\nabla_h \cdot (\roth \boldsymbol{\Psi})\right|_{\infty}$ \\
			\hline
			 6 & 1.7412E-03 &    - & 8.4396E-03 &    - & 8.8774E-05 &    - & 6.3726E-04 &    - & 1.58E-14 \\
			12 & 2.1840E-04 & 3.00 & 1.1937E-03 & 2.82 & 5.8039E-06 & 3.94 & 4.5270E-05 & 3.82 & 3.50E-14  \\
			24 & 2.6929E-05 & 3.02 & 1.5170E-04 & 2.98 & 3.5784E-07 & 4.02 & 3.0642E-06 & 3.88 & 2.80E-14  \\
			48 & 3.3416E-06 & 3.01 & 1.9043E-05 & 2.99 & 2.2289E-08 & 4.00 & 1.9532E-07 & 3.97 & 5.42E-14  \\
			\multicolumn{10}{c}{} \\
			\multicolumn{10}{c}{SPDG $\mathcal{O}(5)$} \\
			$N_h$ & $u_{1,L_1}$ & $\mathcal{O}(L_1)$ & $u_{1,L_{\infty}}$ & $\mathcal{O}(L_{\infty})$ & $\Psi_{1,L_1}$ & $\mathcal{O}(L_1)$ & $\Psi_{1,L_{\infty}}$ & $\mathcal{O}(L_{\infty})$ & $\left|\nabla_h \cdot (\roth \boldsymbol{\Psi})\right|_{\infty}$ \\
			\hline
			 6 & 1.7807E-04 &    - & 1.1339E-03 &    - & 6.4311E-06 &    - & 1.8058E-05 &    - & 4.44E-14 \\
			12 & 9.8436E-06 & 4.18 & 7.9305E-05 & 3.84 & 2.0360E-07 & 4.98 & 5.2722E-07 & 5.10 & 7.26E-14  \\
			24 & 5.9312E-07 & 4.05 & 5.0205E-06 & 3.98 & 6.3782E-09 & 5.00 & 1.6267E-08 & 5.02 & 1.53E-14  \\
			48 & 3.6715E-08 & 4.01 & 3.1475E-07 & 4.00 & 1.9946E-10 & 5.00 & 4.9854E-10 & 5.03 & 3.91E-14  \\
			\multicolumn{10}{c}{} \\
			\multicolumn{10}{c}{SPDG $\mathcal{O}(6)$} \\
			$N_h$ & $u_{1,L_1}$ & $\mathcal{O}(L_1)$ & $u_{1,L_{\infty}}$ & $\mathcal{O}(L_{\infty})$ & $\Psi_{1,L_1}$ & $\mathcal{O}(L_1)$ & $\Psi_{1,L_{\infty}}$ & $\mathcal{O}(L_{\infty})$ & $\left|\nabla_h \cdot (\roth \boldsymbol{\Psi})\right|_{\infty}$ \\
			\hline
			 6 & 1.1663E-05 &    - & 6.3438E-05 &    - & 2.0636E-07 &    - & 1.4971E-06 &    - & 6.88E-16 \\
			12 & 3.7928E-07 & 4.94 & 2.1485E-06 & 4.88 & 3.3500E-09 & 5.94 & 2.6233E-08 & 5.83 & 5.73E-15  \\
			24 & 1.1965E-08 & 4.99 & 6.8688E-08 & 4.97 & 5.1546E-11 & 6.02 & 4.4200E-10 & 5.89 & 2.73E-14  \\
			48 &    3.7165E-10 & 5.01 &    2.1403E-09 & 5.00 &    7.9994E-13 & 6.01 &   7.0461E-12 & 5.97 & 3.71E-14  \\
		\end{tabular}}
	\end{center}
	\label{tab.divcurl_conv}
\end{table}

Finally, we test the div-curl operator by assigning a random vector field $\bpsi(\x)$, generated with the intrinsic function \texttt{RANDOM\_NUMBER} of \texttt{Fortran} programming language. The divergence of the curl of $\bpsi$ is computed by using the structure-preserving operator $\Dop_{k \ell}^{\, \pb}$ given by \eqref{eqn.divOpSPDG} as well as the the non compatible divergence $\tilde{\Dop}_{k \ell}^{\, \pb}$ according to \eqref{eqn.divOpDG}. The results are gathered in Table \ref{tab.divSP} from second up to sixth order of accuracy, proving that the correction terms \eqref{eqn.corrDG} are essential to obtain a structure-preserving high order operator.

\begin{table}[!htbp]  
	\caption{Discrete div-curl errors related to a random vector field $\bpsi$ for different degrees $N=\{1,2,3,4,5\}$. The errors are measured in $L_{\infty}$ norm for the structure-preserving operator and the non compatible operator, that is
	$\left|\tilde{\nabla}_h \cdot (\roth \boldsymbol{\Psi})\right|_{\infty}$ and	
		 $\left|\nabla_h \cdot (\roth \boldsymbol{\Psi})\right|_{\infty}$, respectively.}
	\begin{center} 
	\begin{tabular}{c|cc}	
		$N$ & $\left|\tilde{\nabla}_h \cdot (\roth \boldsymbol{\Psi})\right|_{\infty}$ & $\left|{\nabla}_h \cdot (\roth \boldsymbol{\Psi})\right|_{\infty}$ \\
		\hline
		1 & 6.5524E-01 & 1.4067E-15 \\
		2 & 1.5047E+00 & 2.2986E-14 \\
		3 & 3.1548E+00 & 4.3408E-14 \\
		4 & 4.5372E+00 & 1.7260E-13 \\
		5 & 5.8735E+00 & 2.0706E-12 \\
\end{tabular}
\end{center}
\label{tab.divSP}
\end{table}

%--------- END OF SECTION -------------------------------------------------

%--------- SECTION --------------------------------------------------------
\section{Application to the incompressible Navier-Stokes equations} \label{sec:NSvorticity}
The structure-preserving DG div-curl operators will be tested by considering the incompressible Navier-Stokes equations written in vortex-stream formulation (see \ref{app.vortexNS} for the derivation of the model):
\begin{subequations}
	\begin{align}
		\frac{\partial \boo}{\partial t} - \rot(\uu\times \boo)&=\nu \, \nabla^2 \boo, \label{eqn.vort1} \\
		\rot\rot\bpsi&=\boo, \label{eqn.vort2}\\
		\rot \bpsi &= \uu, \label{eqn.vort3} \\
		\boo&=\boo_0\label{eqn.init}.
	\end{align}
	\label{eqn.vort}
\end{subequations} 

\noindent
The fluid velocity is $\uu=(u,v,w)$, and it is retrieved as the curl of a stream potential function $\bpsi=(\Psi_1,\Psi_2,\Psi_3)$, while the vorticity is $\boo=(\omega_1,\omega_2,\omega_3)=\rot \uu$. We assume a three-dimensional computational domain $\Omega(\x)$ in space, while $t \in \mathds{R}^+$ indicates the time coordinate. The coefficient $\nu$ represents the kinematic viscosity of the fluid, which is assumed to be constant. The PDE system \eqref{eqn.vort} is endowed with three divergence-free constraints, that are given by
\begin{equation}
	\nabla \cdot \boo = 0, \qquad \nabla \cdot \uu = 0, \qquad 	\nabla \cdot \bpsi = 0.
	\label{eqn.divcond}
\end{equation}
The system is then initialized with well-prepared initial condition $\boo_0$ \eqref{eqn.init}, meaning that the involutions \eqref{eqn.divcond} are respected at time $t=0$. Let us notice that the only evolution equation is given by \eqref{eqn.vort1}, whereas the other equations represent compatibility conditions and definitions. Furthermore, recalling that $\nabla \cdot \boo=0$, the vector Laplacian term $\nabla^2 \boo$ can be written as
\begin{eqnarray}
	\nabla^2 \boo &=& \nabla(\nabla \cdot \boo) - \rrot \boo \nonumber \\
	              &=& - \rrot \boo.
	\label{eqn.vecLapl}              
\end{eqnarray}

%--------- END OF SECTION -------------------------------------------------

%--------- SECTION --------------------------------------------------------
\subsection{Numerical method} \label{sec:numscheme}
The numerical discretization of the governing equations \eqref{eqn.vort} consists in a discontinuous Galerkin method in space combined with an implicit-explicit (IMEX) time stepping technique. Particular care must be devoted to satisfy the involutions in space, but also in time so that the structure of the model is mimicked at the discrete level. In other words, the divergence-free constraints \eqref{eqn.divcond} must be respected in the whole computational domain at any time.

A cell-centered discretization is used for the vorticity and the stream potential function, while the velocity vector is defined on a corner-staggered mesh. All variables are approximated in space by piecewise polynomials of arbitrary degree $N$ according to \eqref{eqn.fhat}, thus
\begin{equation}
	\boo^{\, \ib} = \phi_{\ell}(\xxi(\x,C^{\, \ib})) \, \hat{\boo}_{\ell}^{\, \ib}, \quad \boldsymbol{\Psi}^{\, \ib} = \phi_{\ell}(\xxi(\x,C^{\, \ib})) \, \hat{\boldsymbol{\, \bpsi}}_{\ell}^{\, \ib}, \quad \uu^{\, \jb} = \phi_{\ell}(\xxi(\x,D^{\, \jb})) \, \hat{\uu}^{\, \jb}_{\ell}. 
	\label{eqn.DGvar}
\end{equation}
Interpolation of quantities from the centered to the staggered grid are carried out relying on the following $L_2-$projection operator
\begin{eqnarray}
	\Aop_{k\ell}^{\mm}(\boldsymbol{\cdot}_{\ell}^{\mm}) &:=& \frac{|C|}{2^d \, \dx \dy \dz} \, \mathbf{M}_{k \ell}^{-1} \, \iiint_0^1 \phi_k\left( \frac{1-\mm}{4}+\frac{\boldsymbol{\chi}}{2} \right)\phi_\ell\left(\frac{\mm+1}{4}+\frac{\boldsymbol{\chi}}{2}\right) \, d\chi_1 d\chi_2 d\chi_3.
	\label{eqn.L2proj}
\end{eqnarray} 
Due to the self-similarity of the Cartesian mesh, the same interpolation operator can be used to switch from primal to dual mesh and viceversa. For instance, the projections of the vector $\boo^{\, \ib}$ can be computed as
\begin{equation}
	\boo^{\, \jb}=\Aop_{k\ell}^{\mm} \, \hat{\boo}^{\mm}_{\ell}, \quad 	\boo^{\, \ib}=\Aop_{k\ell}^{\tilde{\mm}} \, \hat{\boo}^{\tilde{\mm}}_{\ell}, \qquad \Aop_{k\ell}^{\mm} \, \Aop_{k\ell}^{\tilde{\mm}} = \mathbf{I}_{k \ell},
\end{equation}
where $\mathbf{I}_{k \ell}$ represents the identity matrix, and the multi-index $\tilde{\mm}=-\mm$ refers to the corner neighbors of a dual cell, hence $\boo^{\mm}_{k}$ are corner-staggered values while $\boo^{\tilde{\mm}}_{k}$ are cell-centered quantities. Let us also introduce a new corner-staggered variable $\mathbf{E}_\ell^n$, in analogy to Faraday law in electromagnetism:
\begin{equation}
	\mathbf{E}_{\ell}^{\, \jb,n}=\uu_{k}^{\, \jb,n} \times \Aop_{k\ell}^{\mm} \left( \hat{\boo}_{\ell}^{\mm,n} \right).
\end{equation}

Integration of the vorticity equation \eqref{eqn.vort1} over a cell-centered control volume $C^{\, \ib}$ with the ansatz \eqref{eqn.DGvar}, and multiplication by the inverse of the mass matrix yields the following DG numerical scheme
\begin{equation}
	\boo^{\, \ib,n+1}_{\ell} + \nu \, \dt \, \Cop_{k \ell}^{\, \ib} \left( \Cop_{k \ell}^{\, \jb} \left( \hat{\boo}_{\ell}^{\mm,\tau,n+1} \right) \right) = \boo^{\, \ib,n}_{\ell} + \dt \, \Cop_{k \ell}^{\, \ib} \left( \hat{\mathbf{E}}_{\ell}^{\mm,\tau,n} \right),
	\label{eqn.DGvort}
\end{equation}
where the discrete curls are discretized using the operator \eqref{eqn.curlOpDG}. The viscous terms have been discretized implicitly so that the maximum admissible time step is independent of the parabolic restriction of order $\mathcal{O}(\dx^2)$ and is computed as
\begin{equation}
	\dt = \textnormal{CFL} \min \limits_{\Omega} \left( \frac{\max |u|}{\dx} + \frac{\max |v|}{\dy} + \frac{\max |w|}{\dz} \right)^{-1}.
	\label{eqn.timestep}
\end{equation}

The algorithm for the solution of the system \eqref{eqn.vort} is made of the following steps.

\begin{enumerate}
	\item Initialization of the numerical solution with well-prepared initial data. Firstly, the corner-staggered velocity field $\uu^{\, \jb}$ is assigned, then the vorticity is computed with
	\begin{equation}
		\boo^{\, \ib,0}= \Cop_{k \ell}^{\, i} \left( \hat{\uu}_{\ell}^{\mm,\tau,0} \right),
		\label{eqn.WPIC}
	\end{equation}
	which ensures $\Dop_{k\ell}^{\, \jb,0}(\hat{\boo}_{\ell}^{\mm})=\mathcal{O}(\epsilon)$. To obtain an order of accuracy $\mathcal{O}(N+1)$ for the vorticity, the initial velocity field must be approximated with order $\mathcal{O}(N+2)$.
	
	\item The evolution equation for the vorticity \eqref{eqn.vort1} is solved with the DG scheme \eqref{eqn.DGvort}, hence obtaining the new vorticity $\boo^{\, \ib,n+1}$. Notice that a linear system has to be solved, because of the implicit discretization of the viscous terms. The matrix-free GMRES method \cite{GMRES} is used as iterative solver, which stops after reaching a prescribed tolerance on the residual $|\boo_{\Omega}^{r+1}-\boo_{\Omega}^{r}|\leq \lambda_{\omega}$, with $r$ being the iteration number and  $\lambda_{\omega} = 10^{-10}$ representing the chosen tolerance.
	
	\item The stream function potential has to be retrieved from the compatibility condition \eqref{eqn.vort2}, that involves a so-called double curl problem. Unfortunately, it is known that there exist infinitely many solutions, due to the infinite-dimensional kernel of the curl operator. Therefore, we rely on the regularization strategy forwarded in \cite{DeltaReg2012}, which slightly modifies the original problem by adding a perturbation of the solution $\delta\bpsi$:
	\begin{equation}
		\Cop_{k \ell}^{\, \ib} \left( \Cop_{k \ell}^{\, \jb} \left( \hat{\bpsi}_{\ell}^{\mm,\tau,n+1} \right) \right) + \delta \bpsi^{n+1} = \boo^{n+1}.
		\label{eqn.CurlPsi}
	\end{equation} 
    According to \cite{DeltaReg2012}, we set $\delta = \min (\dx,\dy,\dz)^{N+1}$, and the GMRES solver stops when a tolerance is reached, namely $\lambda_{\Psi} = \min (\dx,\dy,\dz)^{N+2}$. Consequently, the stream function potential $\bpsi^{n+1}$ is determined by the solution of system \eqref{eqn.CurlPsi} with an order of accuracy that is consistent with the DG scheme \eqref{eqn.DGvort}. 
	
	\item Once the stream function is known, the velocity is simply updated with \eqref{eqn.vort3}, hence
	\begin{equation}
		\uu^{\, \jb,n+1} = \Cop_{k \ell}^{\, \ib} \left( \hat{\bpsi}_{\ell}^{\mm,\tau,n+1} \right).
	\end{equation}
	
\end{enumerate}

\begin{theorem} \label{lemma_divvort}
	Assuming well-prepared initial data according to \eqref{eqn.WPIC}, the DG scheme \eqref{eqn.DGvort} preserves the involutions \eqref{eqn.divcond} at the discrete level in space and time.
\end{theorem}
\begin{proof}
Let us consider the first time step of the algorithm, thus we set $n=0$. Applying the discrete structure-preserving divergence operator \eqref{eqn.divOpSPDG} to the evolution equation \eqref{eqn.vort1} leads to
\begin{equation}
	\Dop_{k \ell}^{\, \jb} \left( \hat{\boo}^{\, \mm,1}_{\ell} \right) + \nu \, \dt \, \Dop_{k \ell}^{\, \jb} \left( \Cop_{k \ell}^{\, \ib} \left( \Cop_{k \ell}^{\, \jb} \left( \hat{\boo}_{\ell}^{\mm,\tau,1} \right) \right) \right) = \Dop_{k \ell}^{\, \jb} \left( \hat{\boo}^{\, \mm,0}_{\ell} \right) + \dt \, \Dop_{k \ell}^{\, \jb} \left( \Cop_{k \ell}^{\, \ib} \left( \hat{\mathbf{E}}_{\ell}^{\mm,\tau,0} \right) \right).
	\label{eqn.divomega}
\end{equation}
The assumption of well-prepared initial data \eqref{eqn.WPIC} imposes that $\Dop_{k \ell}^{\, \jb} \left( \hat{\boo}^{\, \mm,0}_{\ell} \right)=\mathcal{O}(\epsilon)$. Similarly, the application of the divergence operator to the curl operators yields an error of $\mathcal{O}(\epsilon)$ that is proven by Theorem \ref{lemma_divcurlDG}, therefore
\begin{equation}
	\Dop_{k \ell}^{\, \jb} \left( \Cop_{k \ell}^{\, \ib} \left( \Cop_{k \ell}^{\, \jb} \left( \hat{\boo}_{\ell}^{\mm,\tau,1} \right) \right) \right) = \mathcal{O}(\epsilon), \qquad \Dop_{k \ell}^{\, \jb} \left( \Cop_{k \ell}^{\, \ib} \left( \hat{\mathbf{E}}_{\ell}^{\mm,\tau,0} \right) \right) = \mathcal{O}(\epsilon).
\end{equation}
Consequently, Equation \eqref{eqn.divomega} simplifies to
\begin{equation}
	\Dop_{k \ell}^{\, \jb} \left( \hat{\boo}^{\, \mm,1}_{\ell} \right) = \mathcal{O}(\epsilon).
\end{equation}
The same procedure can be applied for any successive time step $n>1$.
\end{proof}

\paragraph{Compatible numerical dissipation} The DG scheme \eqref{eqn.DGvort} solves the viscous incompressible Navier-Stokes equations, thus there is no need of embedding additional viscosity in the numerical method. However, the vorticity equation \eqref{eqn.vort1} has the same mathematical form of Faraday law in the MHD equations, where shocks and other strong discontinuities arise from the fluid flow. Therefore, for the sake of completeness, we propose a compatible way to supply the DG scheme \eqref{eqn.DGvort} with numerical dissipation. We rely on the artificial viscosity technique, originally introduced in \cite{Chorin67} for incompressible flows. In order to add a compatible numerical viscosity operator, the vector Laplacian \eqref{eqn.vecLapl} is inserted in the curl term $\rot (\uu \times \boo)$ of Eqn.\eqref{eqn.vort1}, hence obtaining
\begin{equation}
\boo^{\, \ib,n+1}_{\ell} + \nu \, \dt \, \Cop_{k \ell}^{\, \ib} \left( \Cop_{k \ell}^{\, \jb} \left( \hat{\boo}_{\ell}^{\mm,\tau,n+1} \right) \right) = \boo^{\, \ib,n}_{\ell} + \dt \, \Cop_{k \ell}^{\, \ib} \left( \hat{\mathbf{E}}_{\ell}^{\mm,\tau,n} \right) - \frac{h}{Re_h} \, \dt \, \Cop_{k \ell}^{\, \ib} \left( \Cop_{k \ell}^{\, \jb} \left( \hat{\boo}_{\ell}^{\mm,\tau,n} \right) \right),
\label{eqn.DGvortVL} 
\end{equation}
where $h=\max(\dx,\dy,\dz)$ is a characteristic mesh spacing, and $Re_h$ is the mesh Reynolds number that is related to the artificial viscosity which is chosen to be added to the scheme (e.g. $Re_h=10^2$). The above scheme is a compatible discretization of $\boo$ with numerical dissipation, since, taking the curl of \eqref{eqn.DGvortVL}, the structure-preserving property demonstrated by Theorem \ref{lemma_divvort} still holds true. The artificial viscosity can also be added to the physical viscosity, so that an implicit discretization is adopted, hence avoiding the numerical dissipation to affect the time step size.

\subsection{High order extension in time}
The numerical scheme \eqref{eqn.DGvort} is arbitrary accurate in space thanks to the DG discretization, while it exhibits only first order accuracy in time. To enforce a high order time stepping discretization which takes into account both explicit and implicit terms, we employ the class of semi-implicit IMEX schemes originally developed in \cite{BosFil2016}. Following this formalism, the vortex-stream incompressible Navier-Stokes equations \eqref{eqn.vort} are written under the form of an autonomous system, that is
\begin{equation}
	\frac{\partial \boo}{\partial t} = \mathcal{H}\left(\boo(t), \boo(t) \right), \qquad \forall t > t_0, \qquad \textnormal{with} \qquad \boo(t_0)= \boo_0,
\end{equation}
with $\boo_0$ representing the initial condition at time $t_0$. The function $\mathcal{H}$ denotes the spatial approximation of the terms $\rot (\uu \times \boo) - \nu \, \nabla^2 \boo$ in \eqref{eqn.vort1}, that is provided by the DG curl and vector Laplacian operators in \eqref{eqn.DGvort}. The first argument of $\mathcal{H}$ is discretized explicitly and it is labeled with $\boo_E$, while the second argument $\boo_I$ is taken implicitly, so that according to \eqref{eqn.DGvort} the function $\mathcal{H}$ explicitly writes
\begin{equation}
	\mathcal{H}(\boo_E,\boo_I) = \rot (\uu_E \times \boo_E) - \nu \, \rrot \boo_I, \qquad \uu_E=\rot \boo_E.
\end{equation} 

High order in time is reached by means of implicit-explicit (IMEX) Runge-Kutta methods \cite{PR_IMEX}. These multi-step time discretizations are based on a total number of $s$ stages that depend on the required accuracy and stability properties. The time marching starts by setting $\boo_E^n=\boo_I^n=\boo^n$, then the stage fluxes for $i = 1, \ldots, s$ are computed in the following way:
\begin{subequations}
	\begin{align}
		\boo_E^i &= \boo_E^n + \dt \sum \limits_{j=1}^{i-1} \tilde{a}_{ij} k_j, \qquad 2 \leq i \leq s, \label{eq.QE} \\[0.5pt]
		\tilde{\boo}_I^i &= \boo_E^n + \dt \sum \limits_{j=1}^{i-1} a_{ij} k_j, \qquad 2 \leq i \leq s, \label{eq.QI}  \\[0.5pt]
		k_i &= \mathcal{H} \left( \boo_E^i, \tilde{\boo}_I^i + \dt \, a_{ii} \, k_i \right), \qquad 1 \leq i \leq s. \label{eq.k} 
	\end{align}
	\label{eqn.IMEXRK}
\end{subequations}
The coefficients of the explicit and implicit Runge-Kutta method are normally described with the double Butcher tableau:
\begin{equation}
	\begin{array}{c|c}
		\tilde{c} & \tilde{A} \\ \hline & \tilde{b}^\top
	\end{array} \qquad
	\begin{array}{c|c}
		c & A \\ \hline & b^\top
	\end{array},
\label{eqn.butcher}
\end{equation}
with the matrices $(\tilde{A},A) \in \R^{s \times s}$ and the vectors $(\tilde{c},c,\tilde{b},b) \in \R^s$. The tilde symbol refers to the explicit scheme and matrix $\tilde{A}=(\tilde{a}_{ij})$ is a lower triangular matrix with zero elements on the diagonal, while $A=({a}_{ij})$ is a triangular matrix which accounts for the implicit scheme, thus having non-zero elements on the diagonal. Let us notice that the computation of the stage fluxes $k_i$ in \eqref{eq.k} corresponds to the solution of the linear system in the DG scheme \eqref{eqn.DGvort}. We adopt stiffly accurate schemes \cite{BosFil2016,BosPar2021}, thus the new solution $\boo_E^{n+1}=\boo_E^{n+1}=\boo^{n+1}$ is simply given by the last stage of the Runge-Kutta algorithm \eqref{eqn.IMEXRK}, i.e. $\boo^{n+1}=\boo^s$. We report in \ref{app.IMEX} the Butcher tableaux corresponding to the IMEX schemes used to discretize the governing PDE \eqref{eqn.vort} up to third order.

%--------- END OF SECTION -------------------------------------------------

%--------- SECTION --------------------------------------------------------
\section{Numerical results} \label{sec:numtest}
In this section, we solve some numerical test problems for the incompressible Navier-Stokes equations in order to validate the IMEX structure-preserving discontinuous Galerkin schemes presented in this article, which are compactly labeled with SPDG. Unless stated otherwise, the CFL number is set to $\textnormal{CFL}=0.9$ in \eqref{eqn.timestep}, and the numerical dissipation is neglected, thus we set a mesh Reynolds number $Re_h=10^{20}$ in the vector Laplacian term \eqref{eqn.DGvortVL}.

\subsection{Numerical convergence studies: 3D Arnold-Beltrami-Childress flow} \label{ssec:conv}
The convergence properties of the novel SPDG schemes are studied by considering the Arnold-Beltrami-Childress (ABC) flow that was originally introduced in \cite{Arnold65,Childress70} as an interesting class of Beltrami flows. The 3D ABC flow represents a non trivial test, in which all velocity components are intrinsically linked one with each other. For this problem an exact solution of the the three-dimensional incompressible Navier-Stokes equations in a periodic domain is known. The analytical velocity and vorticity fields are given by
\begin{equation}
\uu(\xx,t) = \left\{ \begin{array}{l}
\left(\sin(z)+ \cos(y) \right) e^{-\nu t} \\
\left(\sin(x)+ \cos(z) \right) e^{-\nu t} \\
\left(\sin(y)+ \cos(x) \right) e^{-\nu t} 
\end{array}\right. , \qquad \boo(\xx,t) = \left\{ \begin{array}{l}
\left(\cos(y)+ \sin(z) \right) e^{-\nu t} \\
\left(\cos(z)+ \sin(x) \right) e^{-\nu t} \\
\left(\cos(x)+ \sin(y) \right) e^{-\nu t} 
\end{array}\right. .
\label{eqn.ABC_IC}
\end{equation}
The computational domain is chosen to be the cube $\Omega=[-\pi;\pi]^3$, where periodic boundary conditions are set everywhere. The initial condition is given in terms of the vorticity field according to \eqref{eqn.ABC_IC} by setting $t=0$. The same number $N_h$ of control volumes is used to discretize the domain along each spatial direction, i.e. $N_h=N_x=N_y=N_z$, hence giving rise to a computational grid made of a total number of $N_h^3$ cells. The convergence studies are performed by considering a sequence of successively refined meshes with $N_h=\{8,16,24,32\}$, with the characteristic mesh size $h=\dx=\dy=\dz$.

The first study of convergence is carried out analyzing the $\delta-$regularization technique involved in the solution of the linear system \eqref{eqn.CurlPsi} for the determination of the stream function $\bpsi$. Specifically, after assigning the vorticity field $\boo(\xx,0)$, the linear system for the unknown stream function is solved and the following error is evaluated in $L_{\infty}$ norm:
\begin{equation}
	\varepsilon_{\delta}=\max \limits_{i \in N_e} \left| \rroth \bpsi - \boo \right|_{\infty}.
\end{equation}
Table \ref{tab.delta_conv} reports the convergence rates up to third order using two different choices of the $\delta$ parameter, namely $\delta=h^{N+1}$ and $\delta=h^{N}$. According to the theoretical analysis proposed in \cite{DeltaReg2012}, the formal order of accuracy is obtained with the choice $\delta=h^{N+1}$, while one order of convergence is lost if $\delta=h^{N}$. A tolerance $\tau=h^{N+2}$ has been set in the GMRES solver adopted for the solution of the linear system associated to the stream function, so that the error does not affect the convergence order of the underlying numerical scheme.   

\begin{table}[!htbp]  
	\caption{Numerical convergence results of the $\delta$-regularization method. The errors are measured in $L_\infty$ norm and refer to the quantity $|\rrot \bpsi - \boo|$ at the initial time $t=0$.}
	\begin{center} 
		\renewcommand{\arraystretch}{1.05}
	\begin{tabular}{c|cccc}
			\multicolumn{5}{c}{SPDG $\mathcal{O}(2)$} \\
			$N_h$ & $\varepsilon_{\delta}$ ($\delta=h^{N+1})$ & $\mathcal{O}({L_{\infty}})$ & $\varepsilon_{\delta}$ $(\delta=h^{N})$ & $\mathcal{O}({L_{\infty}})$ \\
			\hline
			 8 & 7.5629E-01 &    - & 8.6743E-01 &    - \\
			16 & 2.6617E-01 & 1.51 & 5.6125E-01 & 0.63 \\
			24 & 1.2809E-01 & 1.80 & 4.1416E-01 & 0.75 \\
			32 & 7.4184E-02 & 1.90 & 3.2791E-01 & 0.81 \\
			\multicolumn{5}{c}{} \\
			\multicolumn{5}{c}{SPDG $\mathcal{O}(3)$} \\
			$N_h$ & $\varepsilon_{\delta}$ $(\delta=h^{N+1})$ & $\mathcal{O}({L_{\infty}})$ & $\varepsilon_{\delta}$ $(\delta=h^{N})$ & $\mathcal{O}({L_{\infty}})$ \\
			\hline
			8  & 8.0180E-01 &    - & 8.7586E-01 &    - \\
			16 & 1.1520E-01 & 2.80 & 2.6754E-01 & 1.71 \\
			24 & 3.5277E-02 & 2.92 & 1.2824E-01 & 1.81 \\
			32 & 1.5073E-02 & 2.96 & 7.4232E-02 & 1.90 \\
	\end{tabular}
\end{center}
\label{tab.delta_conv}
\end{table}

\begin{figure}[!htbp]
	\begin{center}
		\begin{tabular}{cc}
			\includegraphics[width=0.47\textwidth]{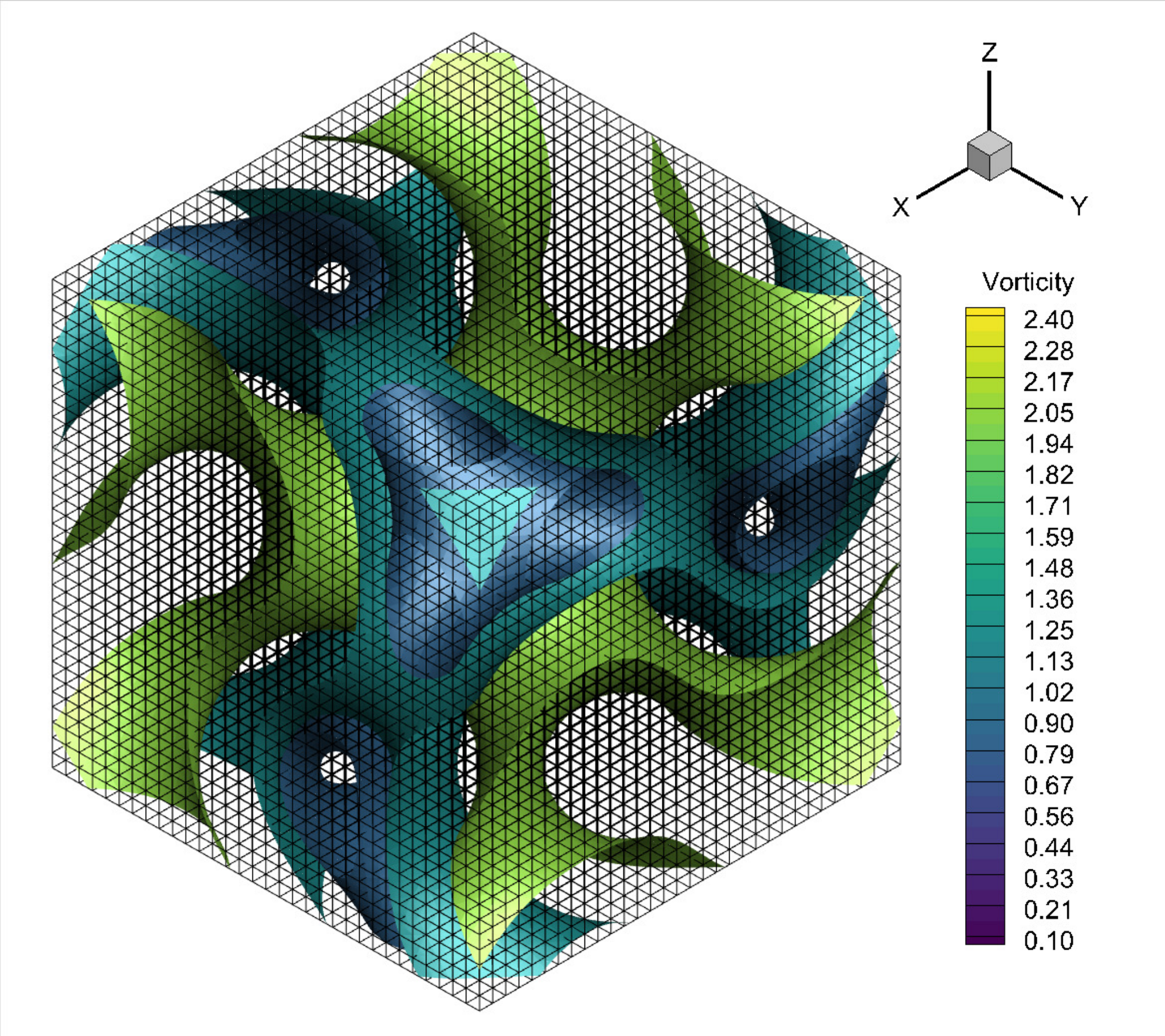}  &          
			\includegraphics[width=0.47\textwidth]{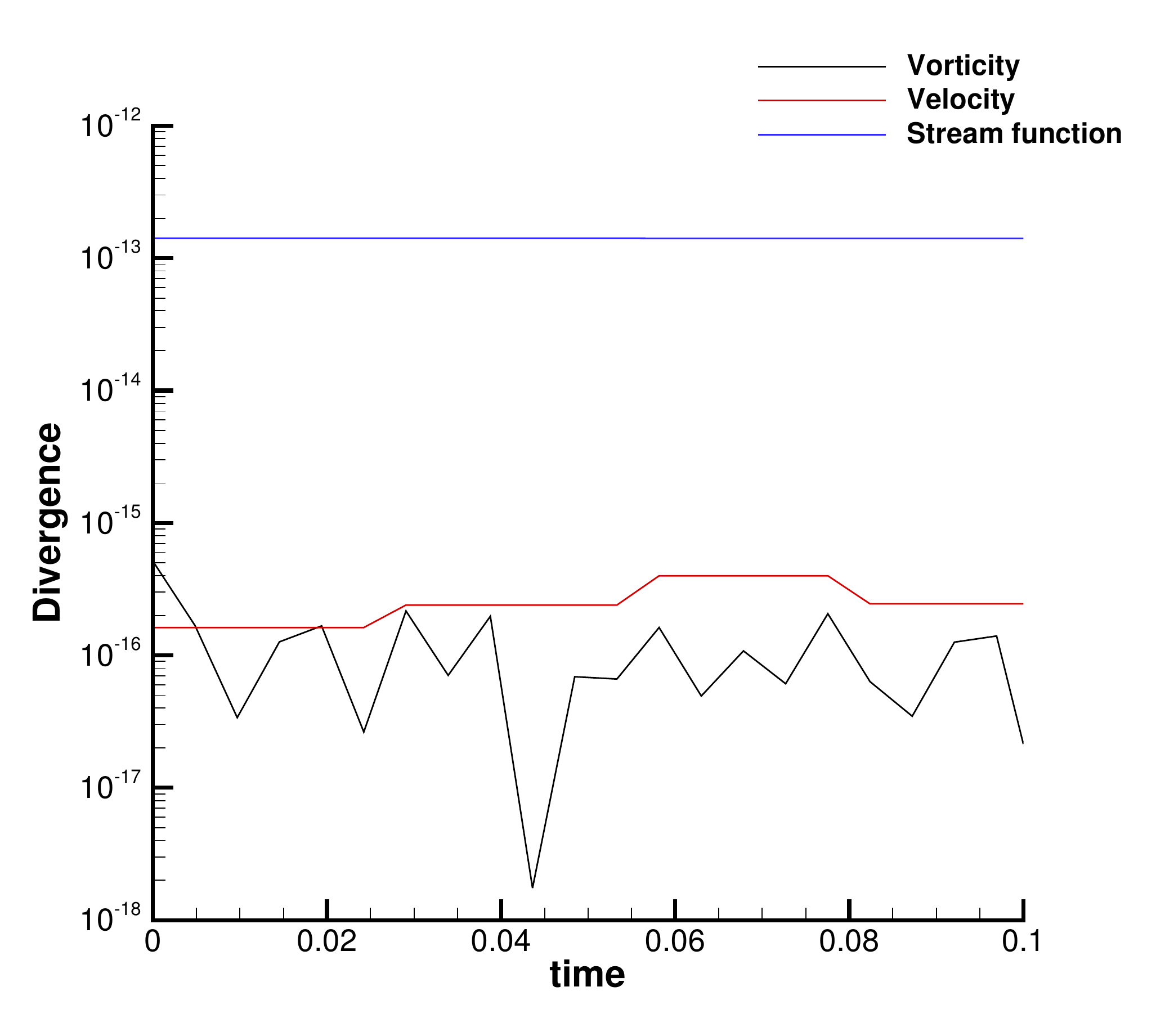}    \\
			\includegraphics[width=0.47\textwidth]{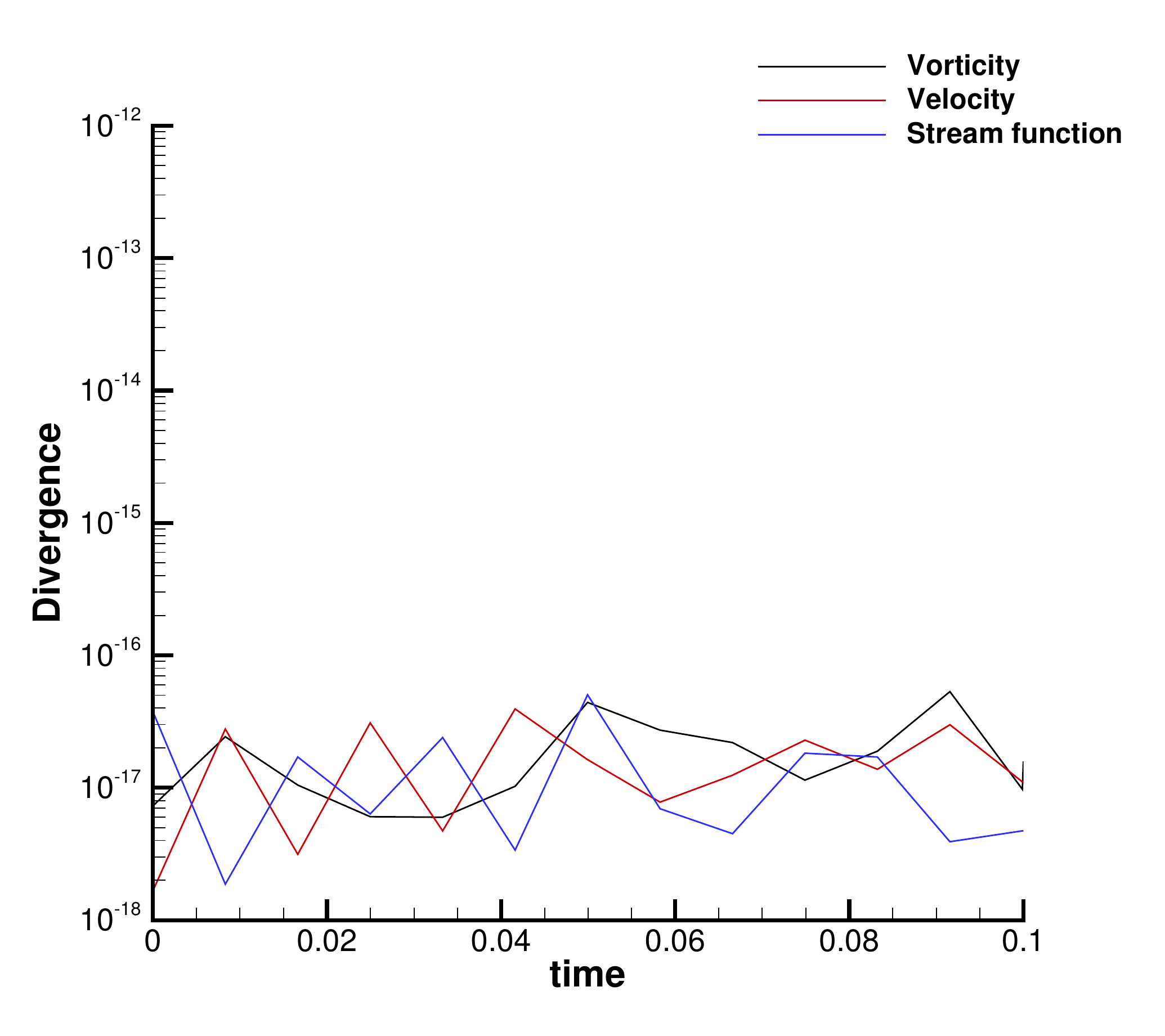}  &          
			\includegraphics[width=0.47\textwidth]{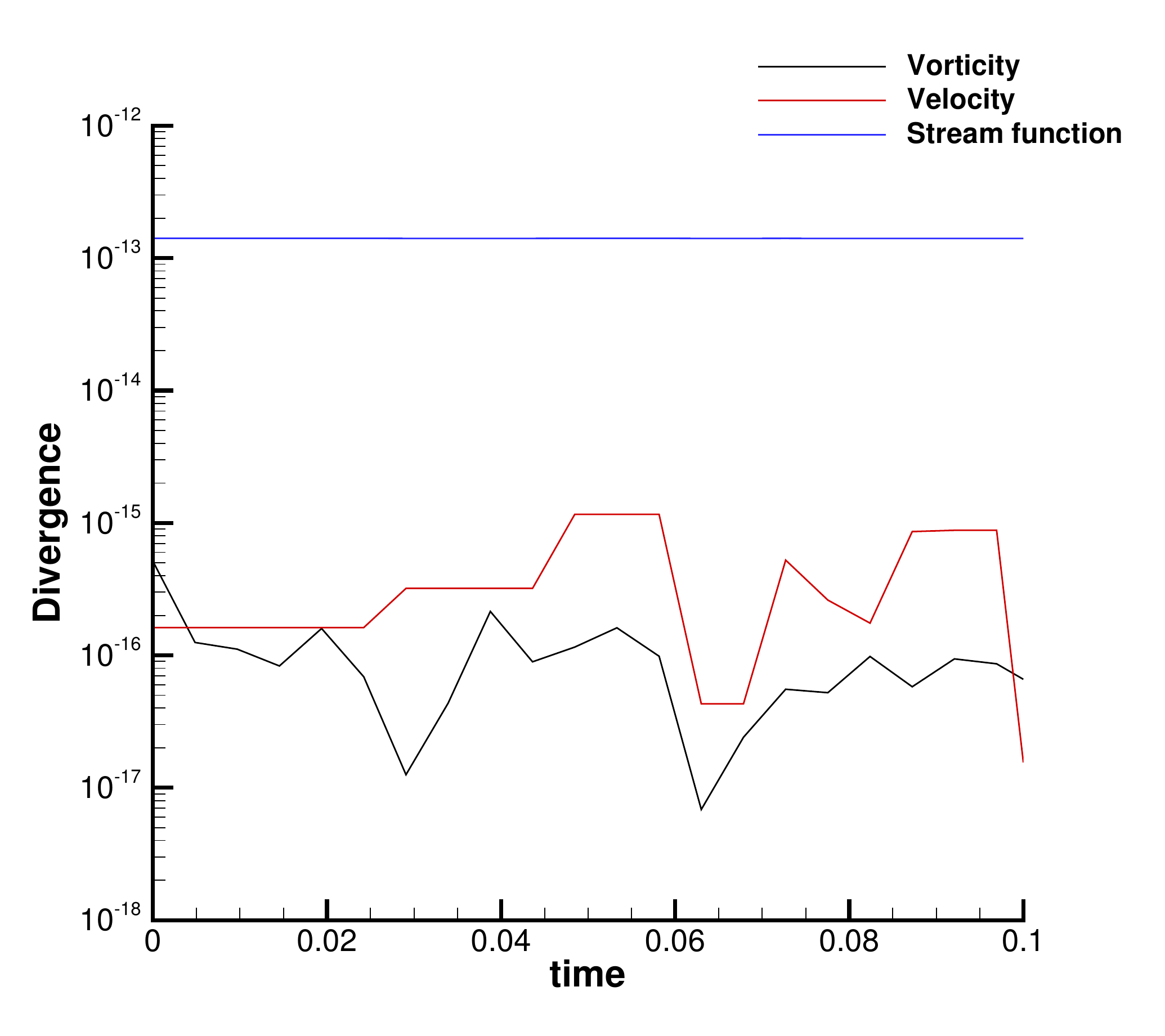}    \\
		\end{tabular}
		\caption{3D Arnold-Beltrami-Childress flow with $\nu=0$. Top: vorticity isosurfaces at levels $[0.8, 1.2, 2.0]$ at time $t=0.1$ (left) and time evolution of the divergence of the vorticity, velocity and stream function fields (right) obtained with third order SPDG schemes. Bottom: time evolution of the divergence errors using the compatible numerical viscosity for second order (left) and third order (right) SPDG schemes with $Re_h=10^2$ and $Re_h=10^4$, respectively.}
		\label{fig.3DABC}
	\end{center}
\end{figure}

The second setup considers the inviscid case, thus the viscosity coefficient is $\nu=0$ and the final time of the simulation is chosen to be $t_f=0.1$. In this case the 3D ABC flow is endowed with a stationary solution, hence allowing to use high order of accuracy in space only. The expected order of convergence is achieved in $L_1$  and $L_{\infty}$ norm for the variables $\omega_1$ and $u$, as confirmed by the analysis shown in Table \ref{tab.inviscid_conv}. Figure \ref{fig.3DABC} shows the vorticity contours at the final time of the simulation as well as the time evolution of the involutions related to vorticity, velocity
and stream function, proving that the SPDG schemes are capable of preserving the divergence-free conditions of the governing equations at the discrete level. To properly check the compatible discretization of the numerical viscosity, the same simulations are also run using the vector Laplacian dissipation \eqref{eqn.DGvortVL} with a mesh Reynolds number of $Re_h=10^2$ and $Re_h=10^4$ for 
$N=1$ and $N=2$, respectively. The results are listed in Table \ref{tab.inviscid_conv_veclapl} confirming that the SPDG schemes can reach third order of accuracy even in the presence of numerical dissipation. Even in this case the divergence-free constraints are respected up to machine precision as depicted in Figure \ref{fig.3DABC}.

\begin{table}[!htbp]  
	\caption{Numerical convergence results of the SPDG schemes \textit{without} numerical viscosity obtained solving the 3D Arnold-Beltrami-Childress flow with $\nu=0$. The errors are measured in $L_1$ and $L_\infty$ norms, and refer to the variables $\omega_1$ and $u_1$ at the final time $t_f=0.1$.}
	\begin{center} 
		\renewcommand{\arraystretch}{1.05}
		\begin{tabular}{c|cccccccc}
			\multicolumn{9}{c}{SPDG $\mathcal{O}(2)$} \\
			$N_h$ & $\omega_{1,L_1}$ & $\mathcal{O}(L_1)$ & $\omega_{1,L_{\infty}}$ & $\mathcal{O}(L_{\infty})$ & $u_{L_1}$ & $\mathcal{O}(L_1)$ & $u_{L_{\infty}}$ & $\mathcal{O}(L_{\infty})$ \\
			\hline
			 8 & 4.6934E+00 &    - & 4.9640E-02 &    - & 7.7848E+01 &    - & 7.5035E-01 &    - \\
			16 & 1.3584E+00 & 1.79 & 2.0301E-02 & 1.29 & 2.6770E+01 & 1.54 & 2.7436E-01 & 1.45 \\
			24 & 7.5742E-01 & 1.44 & 1.2939E-02 & 1.11 & 1.2955E+01 & 1.79 & 1.3286E-01 & 1.79 \\
			32 & 5.2525E-01 & 1.27 & 9.3908E-03 & 1.11 & 7.5845E+00 & 1.86 & 7.7098E-02 & 1.89 \\
			\multicolumn{9}{c}{} \\
			\multicolumn{9}{c}{SPDG $\mathcal{O}(3)$} \\
			$N_h$ & $\omega_{1,L_1}$ & $\mathcal{O}(L_1)$ & $\omega_{1,L_{\infty}}$ & $\mathcal{O}(L_{\infty})$ & $u_{L_1}$ & $\mathcal{O}(L_1)$ & $u_{L_{\infty}}$ & $\mathcal{O}(L_{\infty})$ \\
			\hline
			 8 & 4.5877E-01 &    - & 1.5551E-02 &    - & 6.6853E+01 &    - & 6.7124E-01 &    - \\
			16 & 4.8155E-02 & 3.25 & 9.1384E-04 & 4.09 & 1.1479E+01 & 2.54 & 1.1433E-01 & 2.55 \\
			24 & 1.6543E-02 & 2.64 & 3.3976E-04 & 2.44 & 3.5435E+00 & 2.90 & 3.5258E-02 & 2.90 \\
			32 & 8.2083E-03 & 2.44 & 1.8924E-04 & 2.03 & 1.5104E+00 & 2.96 & 1.5027E-02 & 2.96 \\
		\end{tabular}
	\end{center}
	\label{tab.inviscid_conv}
\end{table}

\begin{table}[!htbp]  
	\caption{Numerical convergence results of the SPDG schemes \textit{with} compatible numerical viscosity obtained solving the 3D Arnold-Beltrami-Childress flow with $\nu=0$. The errors are measured in $L_1$ and $L_\infty$ norms, and refer to the variables $\omega_1$ and $u_1$ at the final time $t_f=0.1$. The mesh Reynolds number is indicated with $Re_h$.}
	\begin{center} 
		\renewcommand{\arraystretch}{1.05}
		\begin{tabular}{c|cccccccc}
			\multicolumn{9}{c}{SPDG $\mathcal{O}(2)$ with $Re_h=10^2$} \\
			$N_h$ & $\omega_{1,L_1}$ & $\mathcal{O}(L_1)$ & $\omega_{1,L_{\infty}}$ & $\mathcal{O}(L_{\infty})$ & $u_{L_1}$ & $\mathcal{O}(L_1)$ & $u_{L_{\infty}}$ & $\mathcal{O}(L_{\infty})$ \\
			\hline
			 8 & 4.7727E+00 &    - & 5.3261E-02 &    - & 7.6864E+01 &    - & 7.4140E-01 &    - \\
			16 & 1.4309E+00 & 1.74 & 1.9787E-02 & 1.43 & 2.6917E+01 & 1.51 & 2.7634E-01 & 1.42 \\
			24 & 8.2275E-01 & 1.36 & 1.2229E-02 & 1.19 & 1.3159E+01 & 1.77 & 1.3535E-01 & 1.76 \\
			32 & 5.7918E-01 & 1.22 & 8.8393E-03 & 1.13 & 7.7607E+00 & 1.84 & 7.9303E-02 & 1.86 \\
			\multicolumn{9}{c}{} \\
			\multicolumn{9}{c}{SPDG $\mathcal{O}(3)$ with $Re_h=10^4$} \\
			$N_h$ & $\omega_{1,L_1}$ & $\mathcal{O}(L_1)$ & $\omega_{1,L_{\infty}}$ & $\mathcal{O}(L_{\infty})$ & $u_{L_1}$ & $\mathcal{O}(L_1)$ & $u_{L_{\infty}}$ & $\mathcal{O}(L_{\infty})$ \\
			\hline
			 8 & 6.8361E-01 &    - & 1.8607E-02 &    - & 6.6784E+01 &    - & 6.5548E-01 &    - \\
			16 & 4.9910E-02 & 3.78 & 1.1244E-04 & 4.05 & 1.1483E+01 & 2.54 & 1.1428E-01 & 2.52 \\
			24 & 1.7778E-02 & 2.55 & 3.6750E-04 & 2.76 & 3.5447E+00 & 2.90 & 3.5265E-02 & 2.90 \\
			32 & 9.5210E-03 & 2.17 & 2.0129E-04 & 2.09 & 1.5105E+00 & 2.97 & 1.5028E-02 & 2.97 \\
		\end{tabular}
	\end{center}
	\label{tab.inviscid_conv_veclapl}
\end{table}

Finally, a fully space-time convergence study is performed by running the 3D Arnold-Beltrami-Childress flow with $\nu=10^{-2}$, so that the solution is no longer stationary and the time discretization actively contributes to the achievement of the formal order of accuracy. We therefore rely on the second and third order IMEX schemes \eqref{eqn.IMEX2}-\eqref{eqn.IMEX3} combined with the SPDG space discretization. No additional numerical viscosity is added ($Re_h=10^{20}$), and the diffusion terms related to the physical viscosity are solved implicitly according to the compatible scheme \eqref{eqn.DGvortVL}. The errors shown in Table \eqref{tab.viscous_conv} demonstrate that the novel IMEX SPDG schemes fulfill the formal accuracy up to third order in space and time.

\begin{table}[!htbp]  
	\caption{Numerical convergence results of the SPDG schemes obtained solving the 3D Arnold-Beltrami-Childress flow with $\nu=10^{-2}$ and no numerical viscosity ($Re_h=10^{20}$). The errors are measured in $L_1$ and $L_\infty$ norms, and refer to the variables $\omega_1$ and $u_1$ at the final time $t_f=0.1$.}
	\begin{center} 
		\renewcommand{\arraystretch}{1.05}
		\begin{tabular}{c|cccccccc}
			\multicolumn{9}{c}{SPDG $\mathcal{O}(2)$} \\
			$N_h$ & $\omega_{1,L_1}$ & $\mathcal{O}(L_1)$ & $\omega_{1,L_{\infty}}$ & $\mathcal{O}(L_{\infty})$ & $u_{L_1}$ & $\mathcal{O}(L_1)$ & $u_{L_{\infty}}$ & $\mathcal{O}(L_{\infty})$ \\
			\hline
			 8 & 4.7014E+00 &    - & 4.9880E-02 &    - & 7.6665E+01 &    - & 7.3955E-01 &    - \\
			16 & 1.3611E+00 & 1.79 & 2.0190E-02 & 1.30 & 2.6676E+01 & 1.52 & 2.7338E-01 & 1.44 \\
			24 & 7.5681E-01 & 1.45 & 1.2981E-02 & 1.09 & 1.2912E+01 & 1.79 & 1.3243E-01 & 1.79 \\
			32 & 5.2456E-01 & 1.27 & 9.4188E-03 & 1.11 & 7.5694E+00 & 1.86 & 7.6946E-02 & 1.89 \\
			\multicolumn{9}{c}{} \\
			\multicolumn{9}{c}{SPDG $\mathcal{O}(3)$} \\
			$N_h$ & $\omega_{1,L_1}$ & $\mathcal{O}(L_1)$ & $\omega_{1,L_{\infty}}$ & $\mathcal{O}(L_{\infty})$ & $u_{L_1}$ & $\mathcal{O}(L_1)$ & $u_{L_{\infty}}$ & $\mathcal{O}(L_{\infty})$ \\
			\hline
			 8 & 6.7526E-01 &    - & 1.8134E-02 &    - & 6.6584E+01 &    - & 6.5349E-01 &    - \\
			16 & 4.8038E-02 & 3.81 & 8.8275E-04 & 4.36 & 1.1373E+01 & 2.55 & 1.1336E-01 & 2.53 \\
			24 & 1.5856E-02 & 2.73 & 3.0215E-04 & 2.64 & 3.5231E+00 & 2.89 & 3.5059E-02 & 2.89 \\
			32 & 7.2446E-03 & 2.72 & 1.4319E-04 & 2.60 & 1.5062E+00 & 2.95 & 1.4987E-02 & 2.95 \\
		\end{tabular}
	\end{center}
	\label{tab.viscous_conv}
\end{table}

\subsection{2D Taylor-Green vortex} \label{ssec:TGV2D}
\begin{figure}[!htbp]
	\begin{center}
		\begin{tabular}{cc}
			\includegraphics[width=0.47\textwidth]{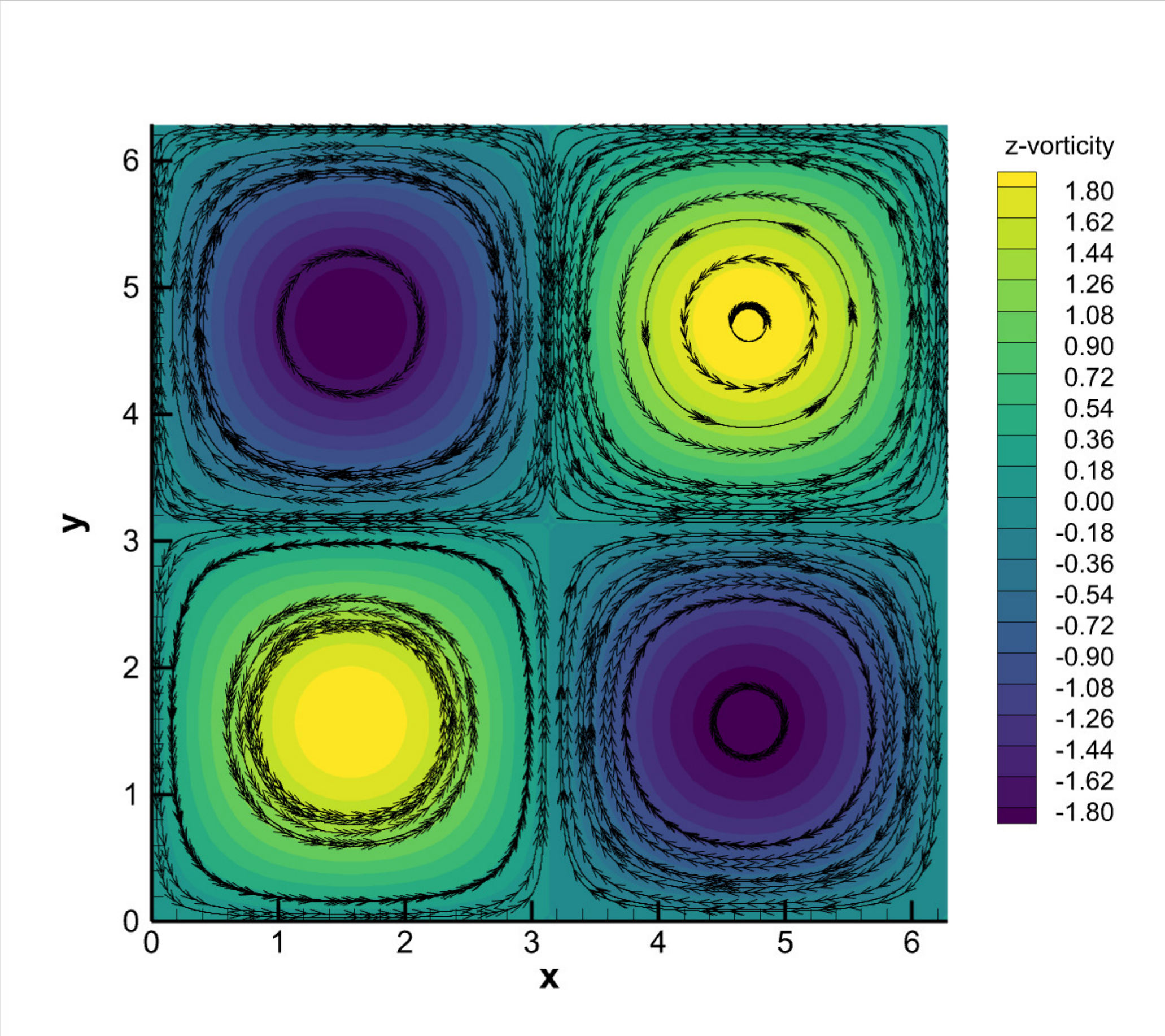}  &          
			\includegraphics[width=0.47\textwidth]{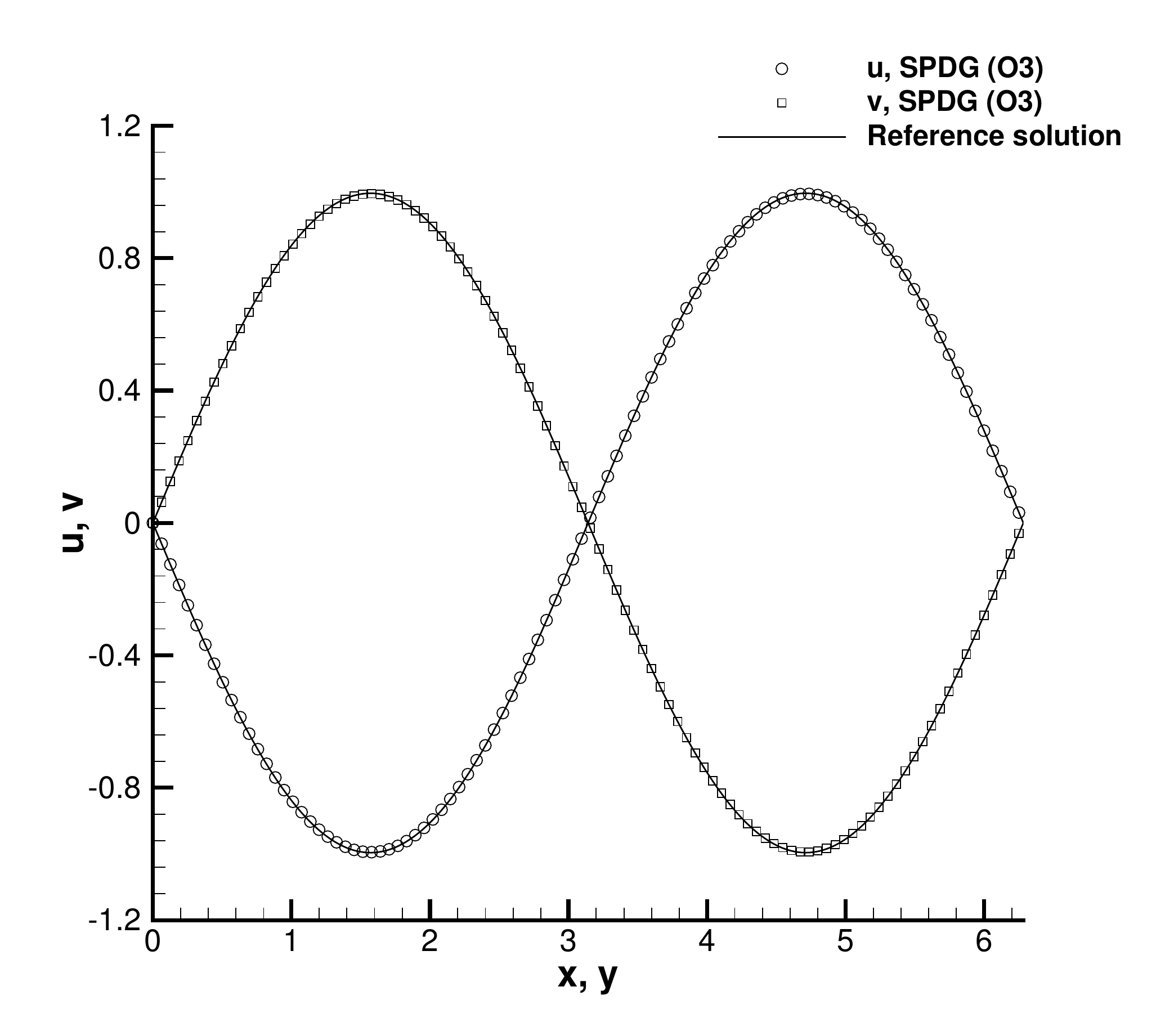}    \\
			\includegraphics[width=0.47\textwidth]{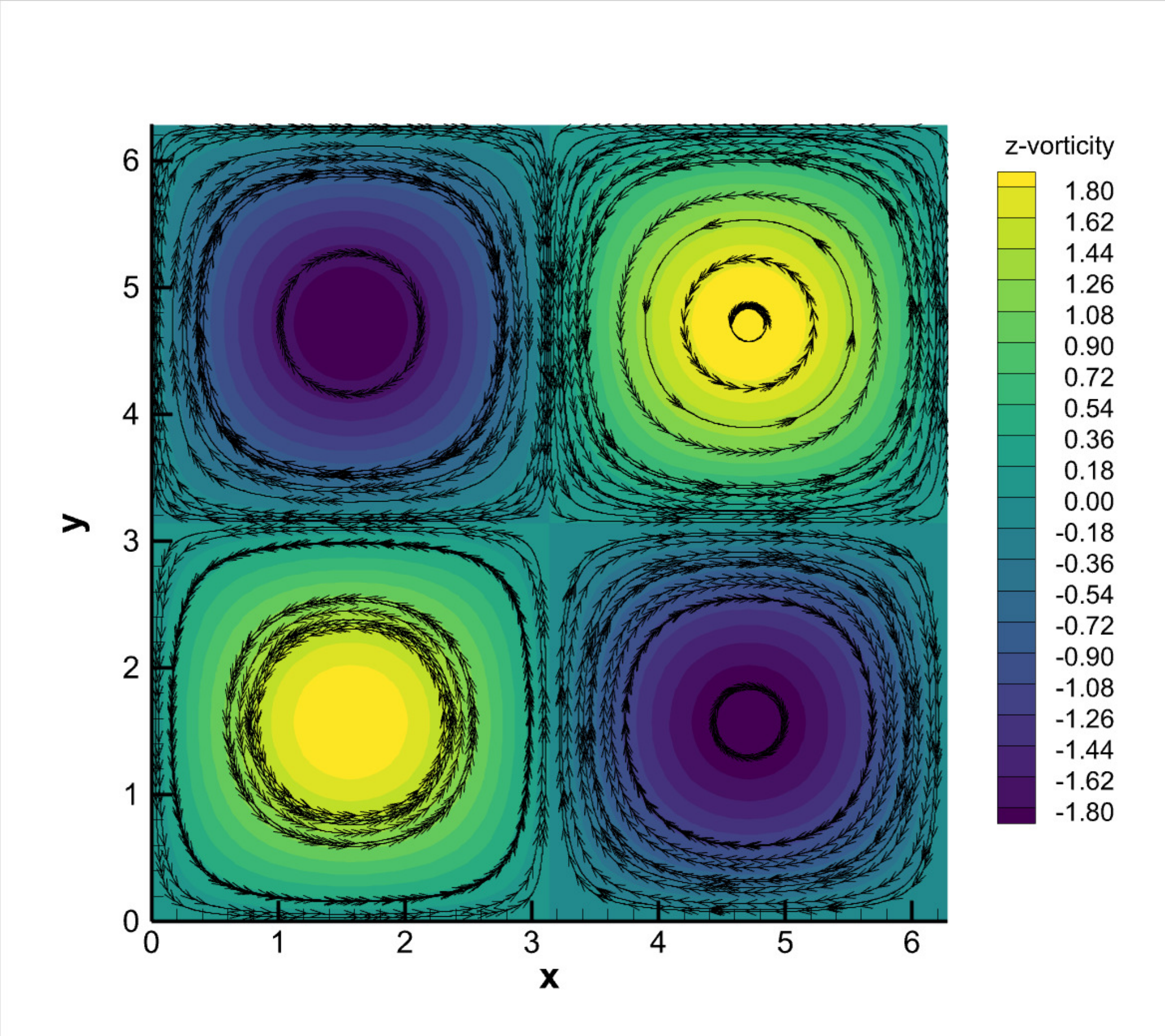}  &          
			\includegraphics[width=0.47\textwidth]{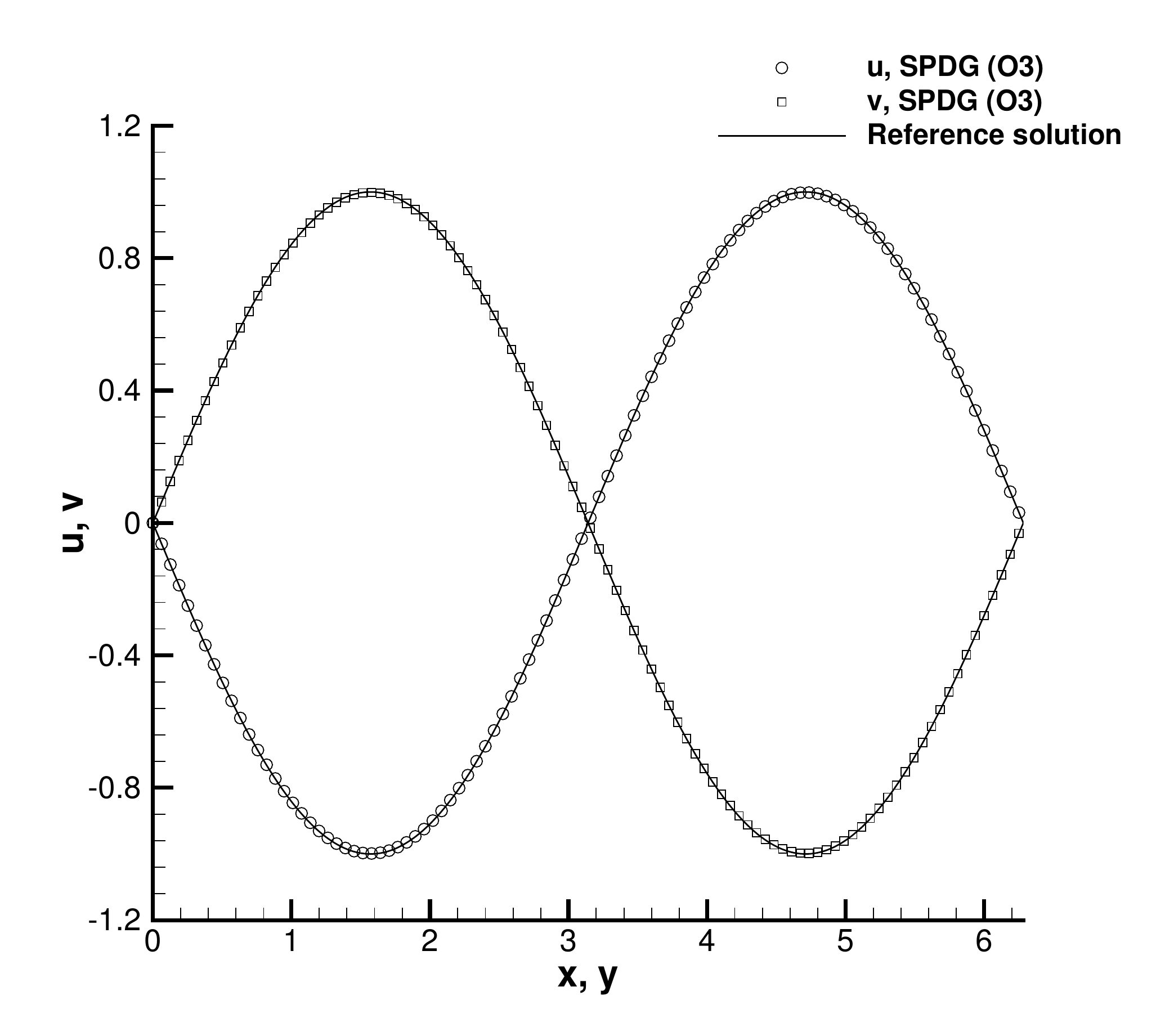}    \\
		\end{tabular}
		\caption{2D Taylor-Green vortex problem at $t_f=0.2$ and comparison with the exact solution of the incompressible Navier-Stokes equations with $\nu=10^{-2}$ (top row) and $\nu=10^{-5}$ (middle row). $z$-vorticity distribution with stream-traces of the velocity field (left) and one-dimensional cuts with 200 equidistant points along the $x-$ and the $y-$axis for velocity components $u$ and $v$ (right).}
		\label{fig.TGV}
	\end{center}
\end{figure}

\begin{figure}[!htbp]
	\begin{center}
		\begin{tabular}{cc}
			\includegraphics[width=0.47\textwidth]{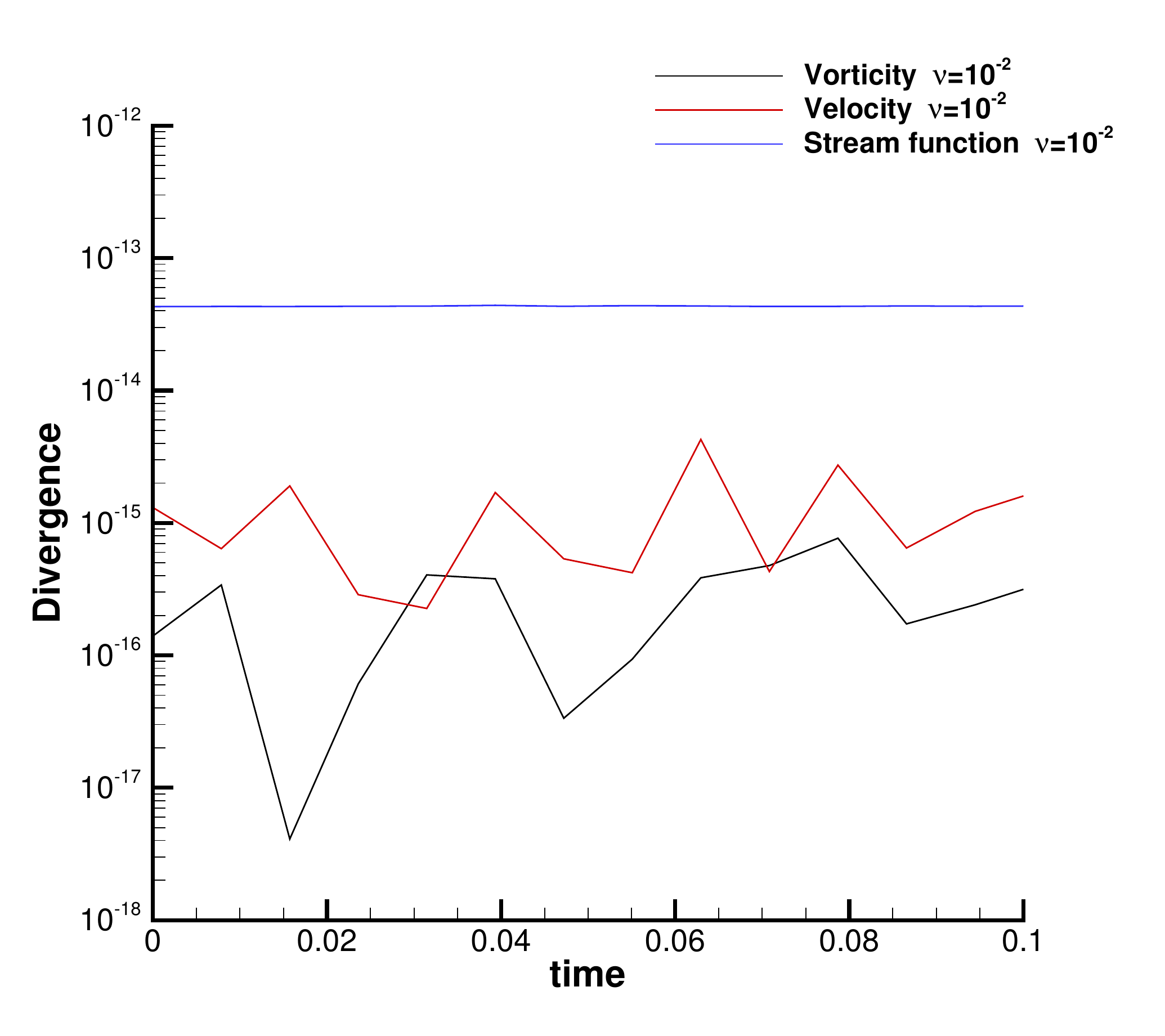}  &          
			\includegraphics[width=0.47\textwidth]{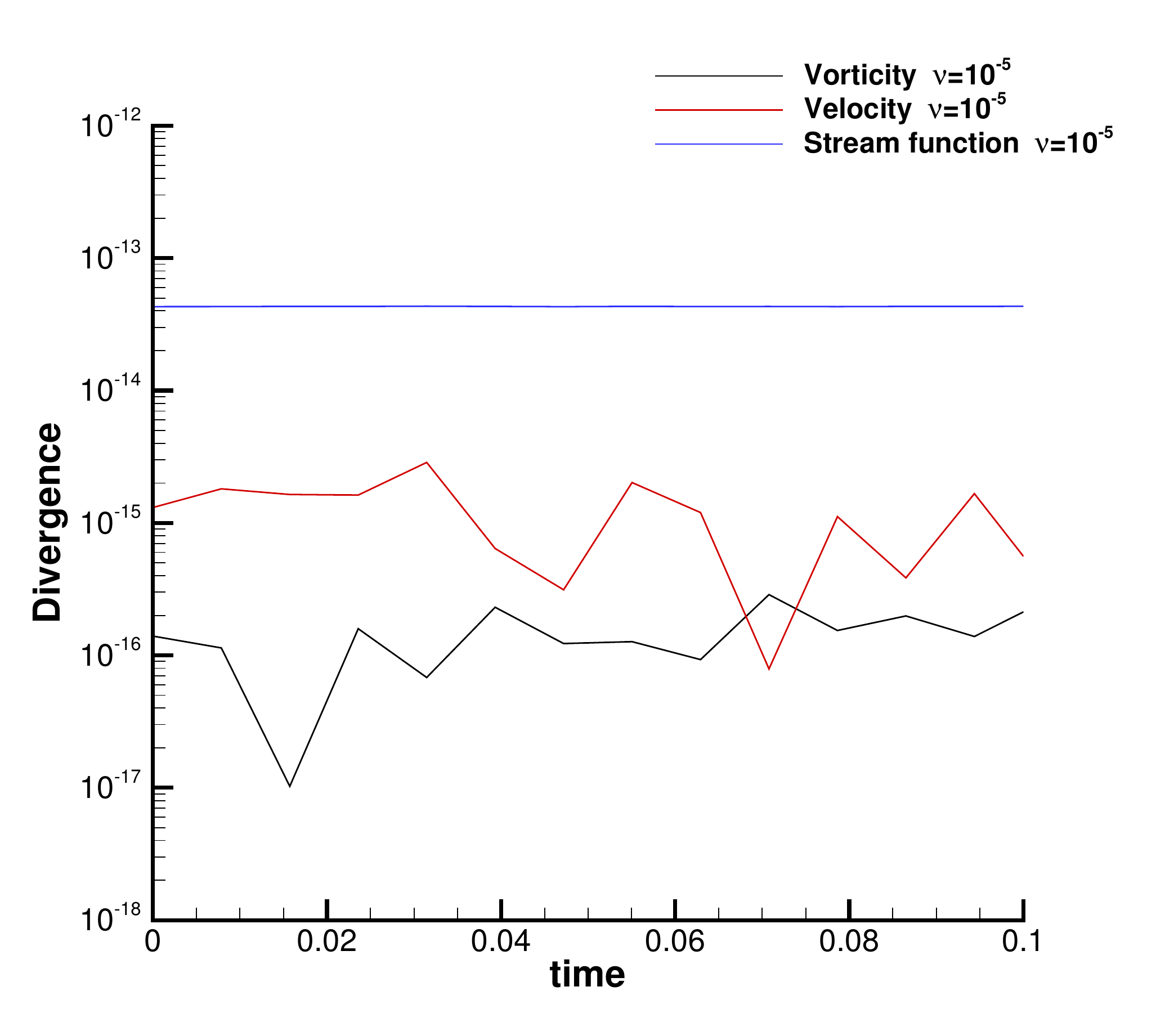}    \\
		\end{tabular}
		\caption{2D Taylor-Green vortex problem. Time evolution of the divergence errors related to vorticity, velocity and stream function with viscosity $\nu=10^{-2}$ (left) and $\nu=10^{-5}$ (right).}
		\label{fig.TGV_divErr}
	\end{center}
\end{figure}

The 2D Taylor-Green vortex represents another rare example for which an exact solution of the incompressible Navier-Stokes equations can be computed, which reads
\begin{equation}
	\uu(\xx,t) = \left\{ \begin{array}{l}
		\sin(x) \cos(y) \, e^{-2 \nu t} \\
		\cos(x) \sin(y) \, e^{-2\nu t} \\
		0 
	\end{array}\right. , \qquad \boo(\xx,t) = \left\{ \begin{array}{l}
		0 \\
		0 \\
		2 \, \sin(x) \sin(y) \, e^{-2 \nu t} 
	\end{array}\right. .
	\label{eqn.TGV_IC}
\end{equation}
This test is run in a three-dimensional periodic domain of size $\Omega=[0;2\pi]^2 \times [0;1]$ with the third order accurate SPDG schemes in space and time. The computational domain is discretized with a rather coarse grid made of $N_x \times N_y \times N_z=48 \times 48 \times 4$ cells. Two different values of the kinematic viscosity are considered, namely $\nu=10^{-2}$ and $\nu=10^{-5}$, and the final time of the simulation is set to $t_f=0.2$. Figure \ref{fig.TGV} depicts the third order numerical solution compared against the analytical solution for the velocity field, exhibiting an excellent matching despite the coarse computational mesh. The stream-traces are smoothly recovered from the SPDG schemes and in accordance with the reference solution. Furthermore, the divergence-free constraints are respected up to machine precision throughout the entire simulation as depicted in Figure \ref{fig.TGV_divErr}, and the structure-preserving property is independent of the physical viscosity thanks to the compatible discretization of the vector Laplacian operator \eqref{eqn.DGvort}.

\subsection{Shear flow test} \label{ssec.dsl}
\begin{figure}[!htbp]
	\begin{center}
		\begin{tabular}{cc}
			\includegraphics[width=0.47\textwidth]{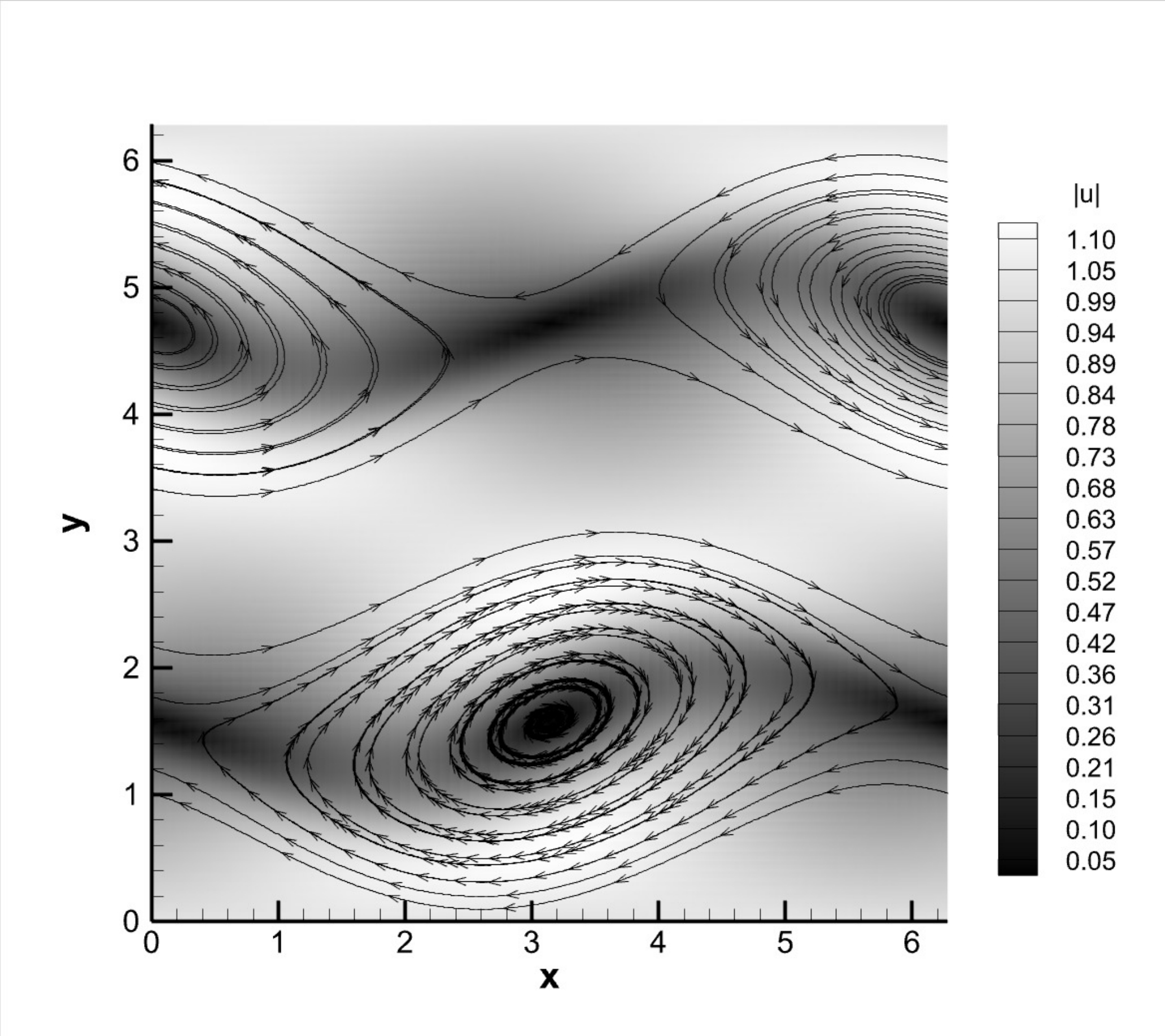}  &          
			\includegraphics[width=0.47\textwidth]{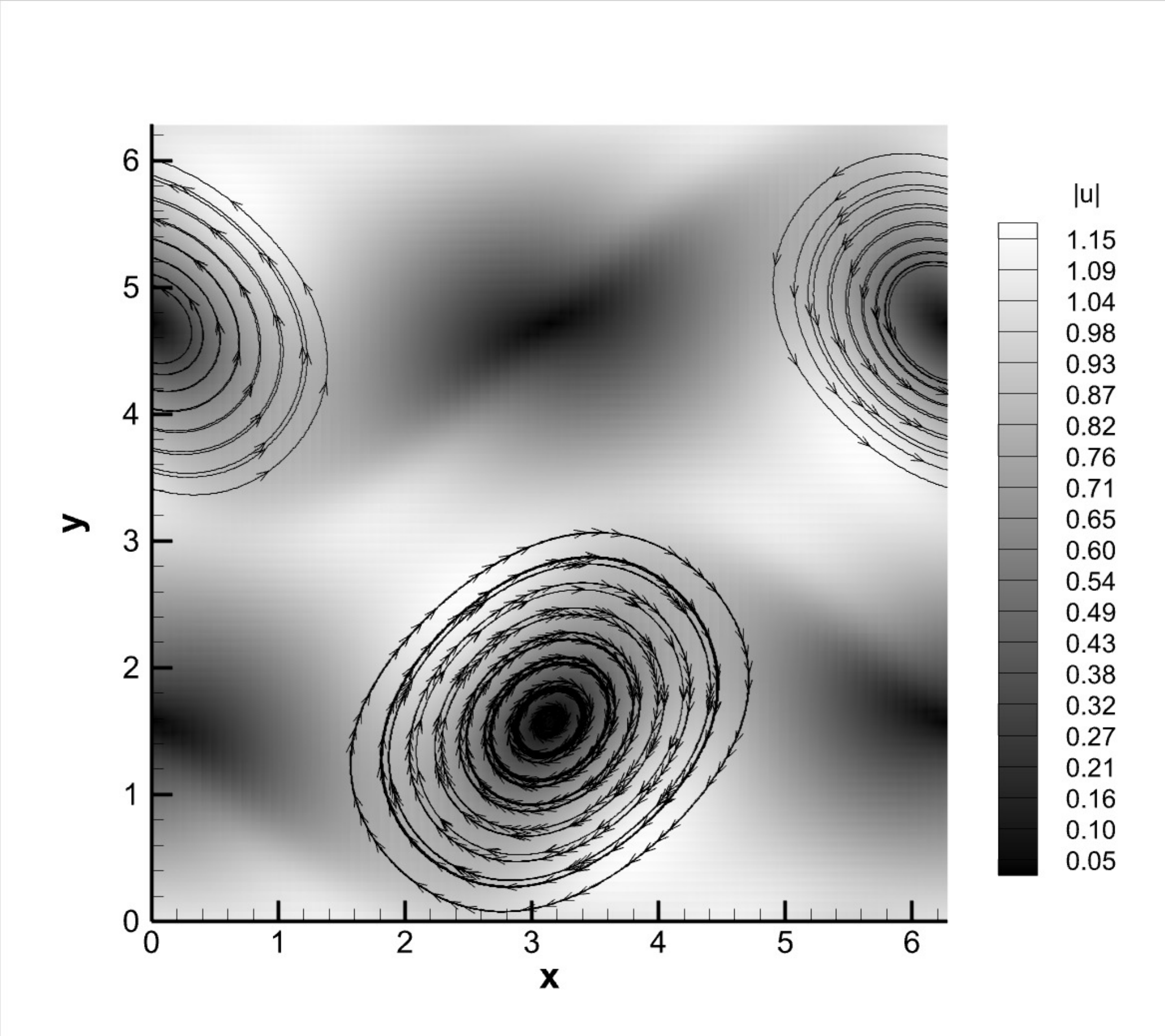}    \\
			\includegraphics[width=0.47\textwidth]{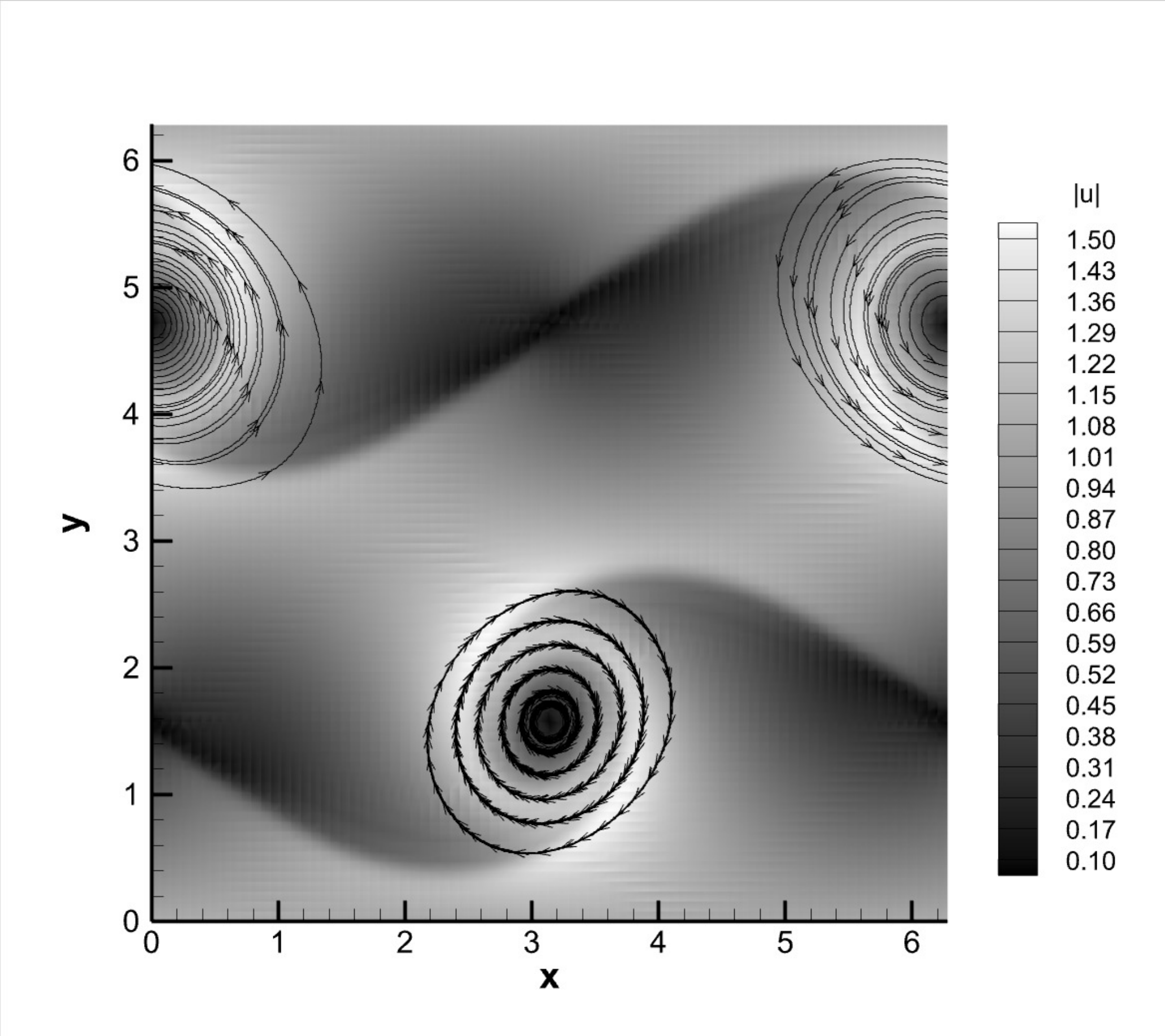}  &          
			\includegraphics[width=0.47\textwidth]{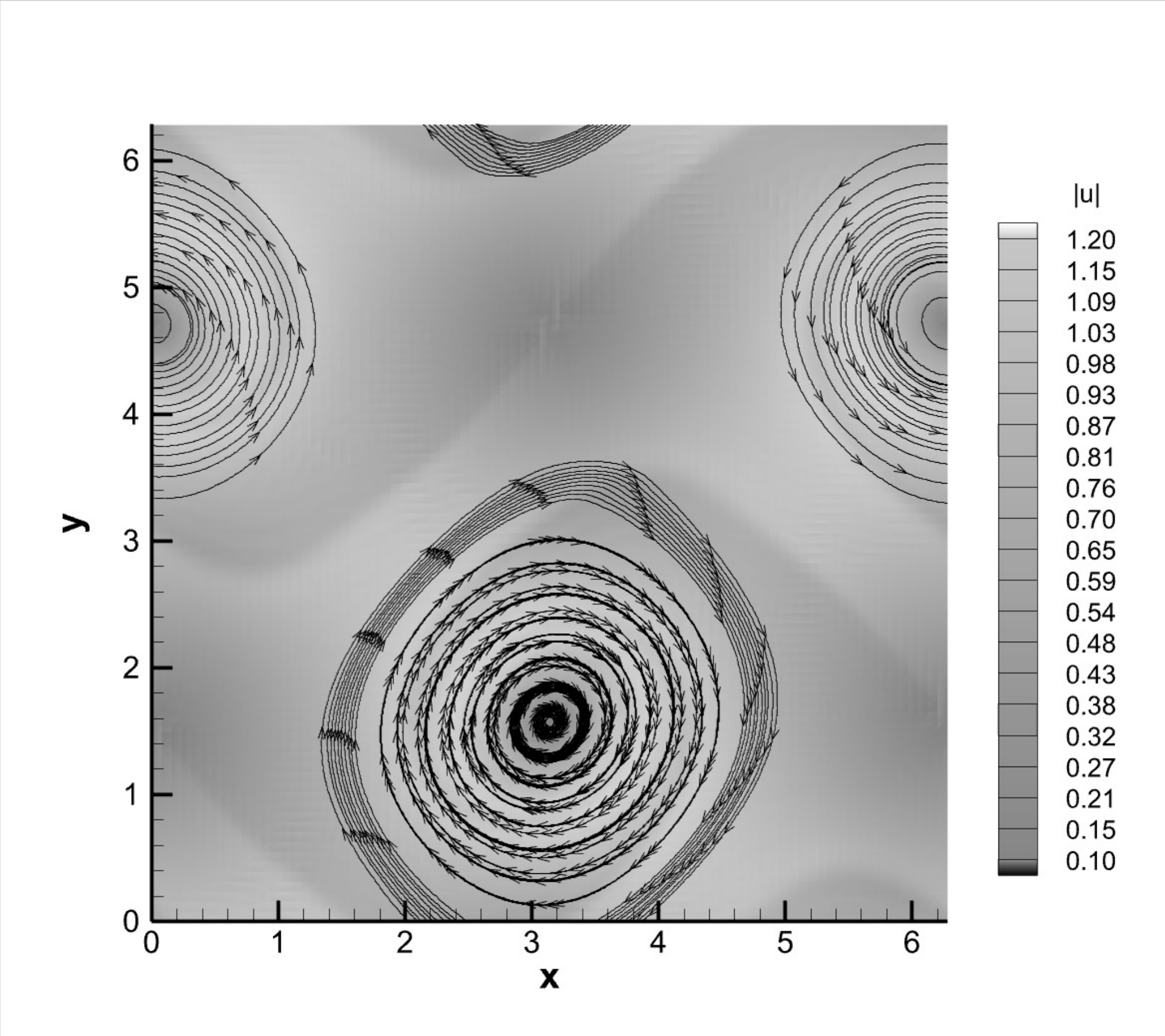}    \\
		\end{tabular}
		\caption{Shear layer test. Velocity magnitude distribution $|\uu|$ with stream-traces of the velocity field at time $t=6$ (left column) and $t=8$ (right column) with viscosity $\nu=10^{-2}$ (top row) and $\nu=10^{-6}$ (bottom row).}
		\label{fig.DSL}
	\end{center}
\end{figure}

The last test case deals with shear flows generated by a horizontal jet that is initially perturbed in a periodic domain. The computational domain is $\Omega=[0,2\pi]^2\times[0,1]$ and it is discretized with a mesh composed of $N_x \times N_y \times N_z=80 \times 80 \times 4$ cubic elements. The initial condition for the vorticity field is imposed according to \cite{QuiShu2011}:
\begin{equation}
	    \omega_1(\xx,t) = \omega_2(\xx,t) = 0, \qquad
		\omega_3(\xx,t) = \left\{ \begin{array}{ll}
		\theta \cos(x) - \frac{1}{\beta} \sech^2((y-\pi/2)/\beta), & y\leq \pi, \\
		\theta \cos(x) + \frac{1}{\beta} \sech^2((3\pi/2-y)/\beta), & y > \pi, 
	\end{array}\right. 
	\label{eqn.TGV_IC}
\end{equation}
with $\theta=0.05$ and $\beta=\pi/15$. The final time of the simulation is $t_f=8$ and this test problem is run using the second order SPDG schemes. Two different viscosity coefficients are chosen, namely an almost inviscid flow with $\nu=10^{-6}$ and a rather viscous fluid with  $\nu=10^{-2}$. Since the viscous terms are discretized implicitly, the maximum admissible time step is not affected by the diffusive eigenvalues but only by the fluid convective speed. Figure \ref{fig.DSL} depicts the stream-traces of the velocity fields at different output times, involving vortexes which arise from the initially perturbed flow layers. The lower the viscous coefficient, the more the vortexes become rolled up, as expected. The divergence errors related to velocity, vorticity and stream function are always maintained at machine accuracy throughout the entire computation, independently on the diffusive terms.

\begin{figure}[!htbp]
	\begin{center}
		\begin{tabular}{cc}
			\includegraphics[width=0.47\textwidth]{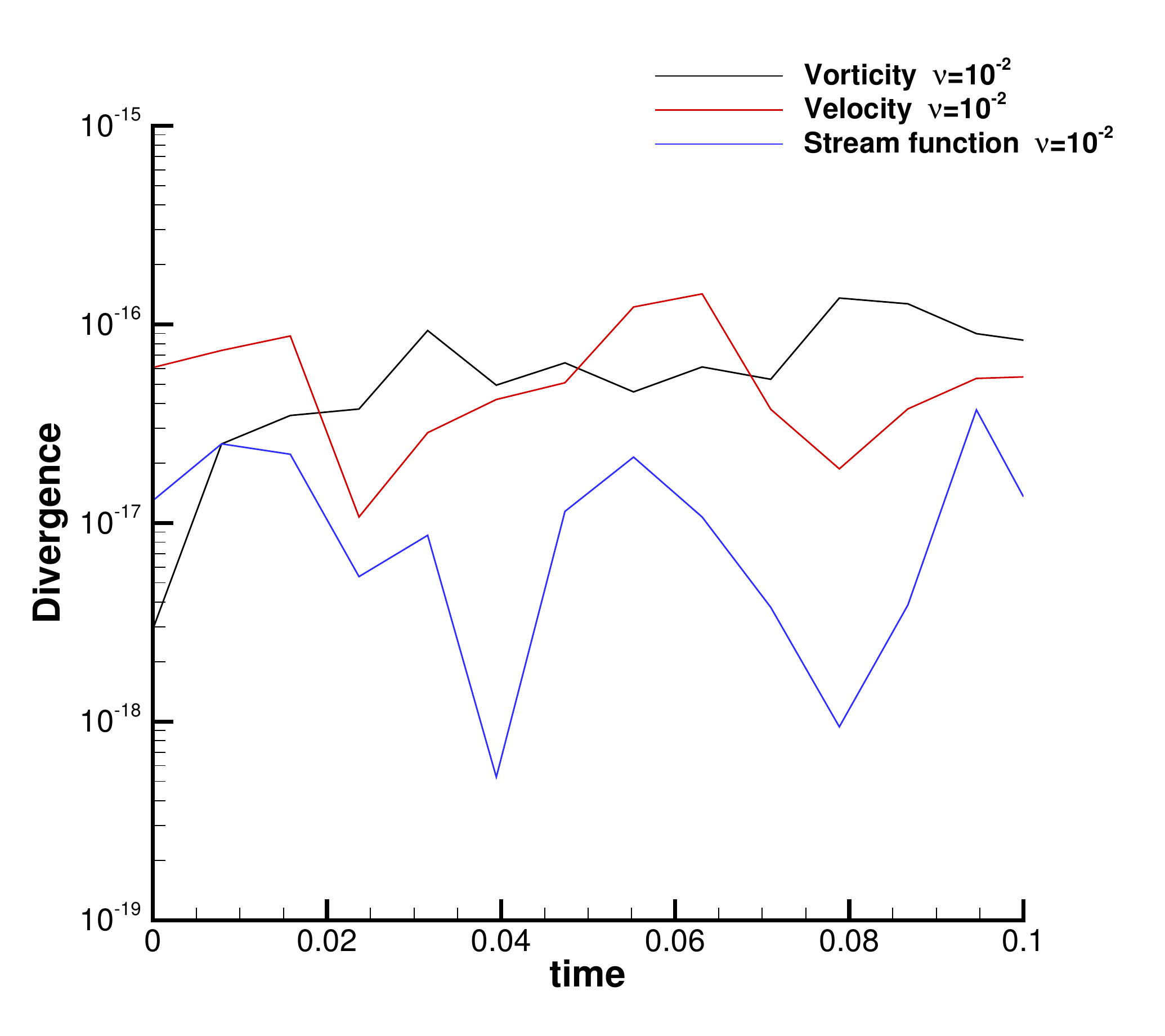}  &          
			\includegraphics[width=0.47\textwidth]{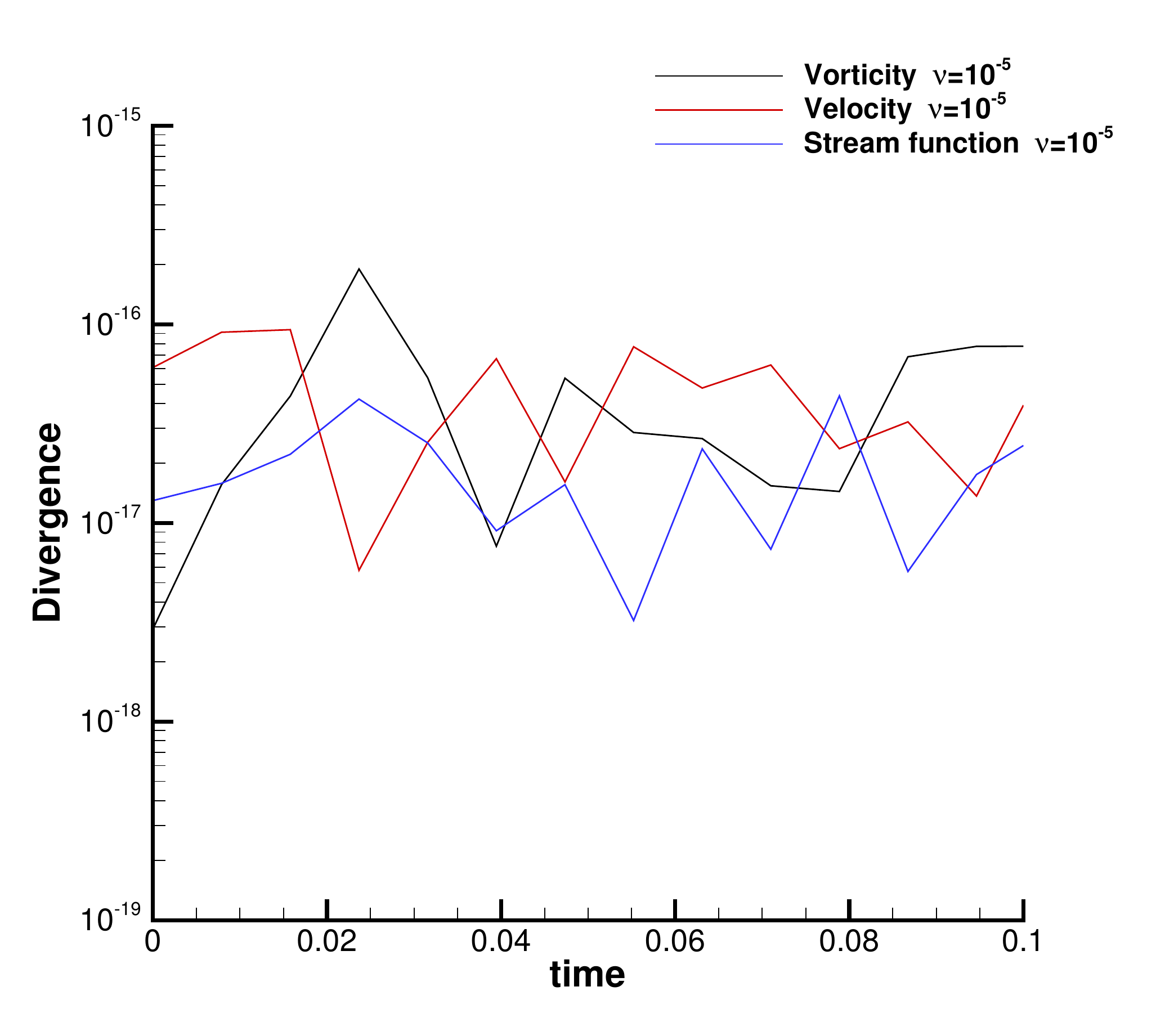}    \\
		\end{tabular}
		\caption{Shear layer test. Time evolution of the divergence errors related to vorticity, velocity and stream function with viscosity $\nu=10^{-2}$ (left) and $\nu=10^{-6}$ (right).}
		\label{fig.DSL_divErr}
	\end{center}
\end{figure}

%--------- END OF SECTION -------------------------------------------------

%--------- SECTION --------------------------------------------------------
\section{Conclusions} \label{sec:conclusion}
A novel structure-preserving scheme has been presented in the discontinuous Galerkin framework. First, the curl operator is written in terms of the divergence operator using the generalized Stokes Theorem, then the discrete divergence is obtained using a corner-staggered discretization on 3D Cartesian grids. The div-curl DG operators collapse by construction to second order structure-preserving finite difference schemes if the DG basis is of degree $N=0$. For higher order schemes, a correction is proposed in order to satisfy Schwarz theorem at the discrete level, demonstrating that the structure-preserving property of the DG div-curl operator is locally and globally satisfied.

The incompressible Navier-Stokes equations written in vortex-stream formulation are used to validate the novel schemes in terms of accuracy and structure-preserving behavior. An implicit-explicit time discretization is adopted to avoid the parabolic time step restriction induced by the diffusion terms, and the novel DG schemes are proven to respect the divergence-free involutions on the vorticity, the velocity and the stream function potential at the discrete level up to machine precision. Particular care has also been devoted to the construction of a compatible numerical viscosity operator, which does not spoil the structure-preserving property of the schemes.

Future work will concern the development of SPDG schemes in the context of plasma flows solving the MHD equations, where the Maxwell-Faraday law can be solved using the same curl operator designed for the vorticity equation. We also plan to investigate curl-free DG operators for applications in continuum mechanics \cite{SIGPR,HOCF2021}.

%--------- END OF SECTION ------------------------------------------------
%
\section*{Acknowledgments}
The authors would like to thank the Italian Ministry of Instruction, University and Research (MIUR) to support this research with funds coming from PRIN Project 2017 (No. 2017KKJP4X entitled “Innovative numerical methods for evolutionary partial differential equations and applications”). 

%--------- END OF SECTION ------------------------------------------------

\newpage
\appendix

\section{Structure-Preserving div-curl finite difference operator} \label{app.SPFD}
Let us introduce a generic vector field $\bar \qq(\xx,t)=(qu,qv,qw)$, which is assumed to be cell-centered without loss of generality. By setting $N=0$, the discrete divergence operator \eqref{eqn.divOpDG} reduces to a central finite difference discretization, since the volume integrals involving $\nabla \phi_k$ vanish and integrals across element boundaries are numerically solved by trapezoidal rule. Therefore, the divergence operator $\nabla_h\cdot (\bar \qq)$ simplifies to
	\begin{equation}
			\nabla_h\cdot(\bar \qq)\ijk = \left(\frac{(qu)_{i+1/2,j,k}-(qu)_{i-1/2,j,k}}{\dx}+\frac{(qv)_{i,j+1/2,k}-(qv)_{i,j-1/2,k}}{\dy}+\frac{(qw)_{i,j,k+1/2}-(qw)_{i,j,k-1/2}}{\dz}\right),
			\label{eqn.divh}
		\end{equation}
with the face-staggered quantities computed by averaging the corresponding values at the centers, i.e.
\begin{equation}
	(qu)_{i+1/2,j,k}=\frac{(qu)_{i+1,j,k}+(qu)_{i,j,k}}{2}, \quad (qv)_{i,j+1/2,k}=\frac{(qv)_{i,j+1,k}+(qv)_{i,j,k}}{2}, \quad (qw)_{i,j,k+1/2}=\frac{(qw)_{i,j,k+1}+(qw)_{i,j,k}}{2}.
	\label{eqn.centraldiff}
\end{equation}
Similarly, the discrete gradient operator \eqref{eqn.gradOpDG} for $N=0$ becomes 
	\begin{equation}
		\nabla(qu)\ijk = \left( \begin{array}{c}
				\frac{qu_{i+1/2,j,k} - qu_{i-1/2,j,k}}{\dx} \\[8pt]
				\frac{qu_{i,j+1/2,k} - qu_{i,j-1/2,k}}{\dy} \\[8pt]
				\frac{qu_{i,j,k+1/2} - qu_{i,j,k-1/2} }{\dz} 
			\end{array} \right),
		\label{eqn.discgrad}
	\end{equation}
and the same holds true for the other components of the vector $\qq$.

We analyze now the discrete curl operator \eqref{eqn.curlOpDG} in the special case $N=0$. After algebraic manipulations, one obtains
	\begin{equation}
			\begin{split}
			\nabla_h\times(\bar \qq)\ijk = &\left(\frac{(qw)_{i,j+1/2,k}-(qw)_{i,j-1/2,k}}{\dy}-\frac{(qv)_{i,j,k+1/2}-(qv)_{i,j,k-1/2}}{\dz},\right.\\
			&\left.\frac{(qu)_{i,j,k+1/2}-(qu)_{i,j,k-1/2}}{\dz}-\frac{(qw)_{i+1/2,j,k}-(qw)_{i-1/2,j,k}}{\dx},\right.\\
			&\left. \frac{(qv)_{i+1/2,j,k}-(qv)_{i-1/2,j,k}}{\dx}-\frac{(qu)_{i,j+1/2,k}-(qu)_{i,j-1/2,k}}{\dy}\right).
		\end{split}
		\label{eqn.curlh}
	\end{equation}

For finite difference schemes, Schwarz theorem is satisfied by the discrete operators \eqref{eqn.divh} and \eqref{eqn.curlh}, thus the following theorem holds true.

\begin{theorem}\label{lemma_divcurl}
		The discrete divergence operator $\nabla_h\cdot(\bar \qq)_{\ijk}$ and curl operator $\nabla_h\times(\bar \qq)\ijk$ satisfy the div-curl property, namely $\nabla_h\cdot\nabla_h\times(\bar \qq)\ijk=0$.
	\end{theorem}
\begin{proof}
		By taking the discrete divergence operator defined in \eqref{eqn.divh} applied to the vector which arises from the curl operator \eqref{eqn.curlh}, each component of the sum is explicitly given as follows:
		\begin{equation}
				\begin{split}
				&\left(\nabla_h\cdot(\bar \qq)\ijk\right)_1 = \frac{(qw)_{i+1/2,j+1/2,k}-(qw)_{i-1/2,j+1/2,k}}{\dx\dy}-\frac{(qw)_{i+1/2,j-1/2,k} -(qw)_{i-1/2,j-1/2,k}}{\dx\dy}\\
				&-\frac{(qv)_{i+1/2,j,k+1/2}-(qv)_{i-1/2,j,k+1/2}}{\dx\dz}+\frac{(qv)_{i+1/2,j,k-1/2}-(qv)_{i-1/2,j,k-1/2}}{\dx\dz},
					\end{split}
			\label{divcurlx}
			\end{equation}
		\begin{equation}
			\begin{split}
					&\left(\nabla_h\cdot(\bar \qq)\ijk\right)_2 = \frac{(qu)_{i,j+1/2,k+1/2}-(qu)_{i,j-1/2,k+1/2}}{\dz\dy}-\frac{(qu)_{i,j+1/2,k-1/2} -(qu)_{i,j-1/2,k-1/2}}{\dz\dy}\\
					&-\frac{(qw)_{i+1/2,j+1/2,k}-(qw)_{i+1/2,j-1/2,k}}{\dx\dy}+\frac{(qw)_{i-1/2,j+1/2,k}-(qw)_{i-1/2,j-1/2,k}}{\dx\dy},
				\end{split}
			\label{divcurly}
		\end{equation}
		\begin{equation}
			\begin{split}
					&\left(\nabla_h\cdot(\bar \qq)\ijk\right)_3 = \frac{(qv)_{i+1/2,j,k+1/2}-(qv)_{i+1/2,j,k-1/2}}{\dz\dx}-\frac{(qv)_{i-1/2,j,k+1/2} -(qu)_{i-1/2,j,k-1/2}}{\dz\dx}\\
					&-\frac{(qu)_{i,j+1/2,k+1/2}-(qu)_{i,j+1/2,k-1/2}}{\dz\dy}+\frac{(qu)_{i,j-1/2,k+1/2}-(qu)_{i,j-1/2,k-1/2}}{\dz\dy}.
				\end{split}
			\label{divcurlz}
		\end{equation}
	Now, by summing up the above defined three quantities one gets
		\begin{equation}
				\left(\nabla_h\cdot(\bar \qq)\ijk\right)_1 +\left(\nabla_h\cdot(\bar \qq)\ijk\right)_2 +\left(\nabla_h\cdot(\bar \qq)\ijk\right)_3=0,
			\label{divcurl}
		\end{equation}
	which corresponds to the divergence free condition.
	\end{proof}

\section{Derivation of the vortex-stream formulation of the incompressible Navier-Stokes equations} \label{app.vortexNS}
The vortex-stream formulation \eqref{eqn.vort} can be obtained starting from the incompressible Navier-Stokes model which reads as
\begin{subequations}
		\begin{align}
				\frac{\partial \uu}{\partial t} + \nabla \cdot \mathbf{F}_c + \nabla p - \nu \, \nabla^2 \uu &=\mathbf{0}.\label{eqn.mom} \\
				\nabla \cdot \uu &=0. \label{eqn.cont}
			\end{align}
		\label{eqn.ins}
	\end{subequations} 
The physical domain is $\Omega \in \mathds{R}^d$ defined in space dimension $d \in \{ 1, 2, 3 \}$ with the vector of spatial coordinates $\xx=(x,y,z) \in \Omega$ and time variable $t \in \mathds{R}^+$. The vector $\uu=(u,v,w)$ is the velocity vector with components along each spatial direction and $p=P/\rho$ is the normalized fluid pressure, $P$ being the physical pressure and $\rho$ denoting the constant fluid density. The kinematic viscosity coefficient is computed by $\nu=\mu/\rho$, and $\mu$ is the dynamic viscosity of the fluid. The flux tensor of the nonlinear convective terms explicitly writes
\begin{equation}
		\mathbf{F}_c:=\uu \otimes \uu = \left[ \begin{array}{ccc}
				uu & uv & uw \\ vu & vv & vw \\ wu & wv & ww
			\end{array} \right].
	\end{equation}
The vorticity $\boo$ is defined by
\begin{equation}\label{vort}
	\rot \uu=\boo=(\oo_x,\oo_y,\oo_z)=\left(\frac{\partial w}{\partial y}-\frac{\partial v}{\partial z},\frac{\partial u}{\partial z}-\frac{\partial w}{\partial x},\frac{\partial v}{\partial x}-\frac{\partial u}{\partial y}\right),
\end{equation}
and an evolution equation for the vorticity is obtained by taking the curl of the momentum equation \eqref{eqn.mom}, that is
\begin{equation}\label{vorticty0}
		\rot \frac{\partial \uu}{\partial t} + \rot (\uu\nabla \cdot \uu) + \rot\nabla p - \nu \, \rot\nabla^2 \uu =\mathbf{0}.
	\end{equation}
Recalling that $\rot\nabla p=0$ and since the inertia term can be rewritten as $\uu\nabla \cdot \uu=\frac{1}{2}\nabla(\uu\cdot\uu)-\uu\times(\rot\uu)$ we get
\begin{equation}
		\frac{\partial \boo}{\partial t} + \uu\cdot\nabla\boo-\boo\cdot\uu=\nu \,\nabla^2 \boo,
	\end{equation}
since by construction $\rot\nabla(\uu\cdot\uu)=0$. The above equation is then coupled with the incompressibility constraint $\nabla \cdot \uu =0$. 

By means of algebraic arguments, the vorticity equation \eqref{vorticty0} can be reformulated as
\begin{equation}\label{vorticty}
		\frac{\partial \boo}{\partial t} - \rot(\uu\times \boo)=\nu \, \nabla^2 \boo,
	\end{equation}
which corresponds to Equation \eqref{eqn.vort1}. 
Once the vorticity is known, since it handles the relation $\rot\boo=\rot\rot \uu=\nabla(\nabla\cdot\uu)-\nabla^2\uu$, the velocity field can be recovered at all times by solving the following Poisson equation
\begin{equation}\label{poisson}
		\nabla^2\uu=-\rot\boo.
	\end{equation}
However, this form is not suitable for devising divergence-free numerical schemes because it is difficult to obey this involution in the Laplacian operator. In alternative, since the velocity field is divergence-free, one can introduce the stream function potential $\bpsi$ such that $\nabla\cdot\bpsi=0$ and then solve the double curl equation
\begin{equation}\label{poisson}
		\rot\rot\bpsi=\boo.
	\end{equation}
Successively the velocity field can be recovered from from the relation $\uu=\rot\bpsi$, so that it is easier to maintain a discrete divergence-free velocity field using a structure-preserving div-curl operator. This is the reason why we solve the vortex-stream incompressible Navier-Stokes equations written in the form \eqref{eqn.vort}.

%=========================================================================
\section{IMEX schemes} \label{app.IMEX}
%=========================================================================
The Butcher tableau \eqref{eqn.butcher} for the IMEX schemes used in this work are reported hereafter. They have been derived in \cite{PR_IMEX,PR_IMEXHO} and each IMEX scheme is described with a triplet $(s,\tilde{s},p)$ which characterizes the number $s$ of stages of the implicit method, the number $\tilde{s}$ of stages of the explicit method and the order $p$ of the resulting scheme. The acronym SA stands for Stiffly Accurate, while DIRK refers to Diagonally Implicit Runge-Kutta schemes.

\begin{itemize}
	\item SP(1,1,1)
	
	\begin{equation}
		\begin{array}{c|c}
			0 & 0 \\ \hline & 1
		\end{array} \qquad
		\begin{array}{c|c}
			1 & 1 \\ \hline & 1
		\end{array}
		\label{eqn.IMEX1}
	\end{equation}
	
	\item LSDIRK2(2,2,2) \hspace{0.2cm} $\gamma=1-1/\sqrt{2}$, \hspace{0.05cm} $\beta=1/(2\gamma)$
	
	\begin{equation}
		\begin{array}{c|cc}
			0 & 0 & 0 \\ \beta & \beta & 0 \\ \hline & 1-\gamma & \gamma
		\end{array} \qquad
		\begin{array}{c|cc}
			\gamma & \gamma & 0 \\ 1 & 1-\gamma & \gamma \\ \hline & 1-\gamma & \gamma
		\end{array}
		\label{eqn.IMEX2}
	\end{equation}
	
	\item SA DIRK (3,4,3) \hspace{0.2cm} $\gamma=0.435866$
	
	\begin{equation}
		\begin{array}{c|cccc}
			0 & 0 & 0 & 0 & 0 \\ \gamma & \gamma & 0 & 0 & 0 \\ 0.717933 & 1.437745 & -0.719812 & 0 & 0 \\ 1 & 0.916993 & 1/2 & -0.416993 & 0 \\ \hline  & 0 & 1.208496 & -0.644363 & \gamma
		\end{array} \qquad
		\begin{array}{c|cccc}
			\gamma & \gamma & 0 & 0 & 0 \\ \gamma & 0 & \gamma & 0 & 0 \\ 0.717933 & 0 & 0.282066 & \gamma & 0  \\ 1 & 0 & 1.208496 & -0.644363 & \gamma \\ \hline  & 0 & 1.208496 & -0.644363 & \gamma
		\end{array}
		\label{eqn.IMEX3}
	\end{equation}

\end{itemize} 

%--------- BIBLIOGRAPHY ------------------------------------------------

\bibliography{biblio}
\bibliographystyle{unsrt}

\end{document}